\documentstyle[10pt, amscd]{amsart}
\renewcommand{\epsilon}{\varepsilon}
\newcommand{\newsection}[1]
{\subsection{#1}\setcounter{theorem}{0} \setcounter{equation}{0}
\par\noindent}

\newtheorem{theorem}{Theorem}

\newtheorem{lemma}[theorem]{Lemma}
\newtheorem{corr}[theorem]{Corollary}

\newtheorem{proposition}[theorem]{Proposition}
\newtheorem{deff}[theorem]{Definition}

\newcommand{\bth}{\begin{theorem}}
\newcommand{\ble}{\begin{lemma}}
\newcommand{\bcor}{\begin{corr}}
\newcommand{\bdeff}{\begin{deff}}
\newcommand{\bprop}{\begin{proposition}}
\newcommand{\eth}{\end{theorem}}
\newcommand{\ele}{\end{lemma}}
\newcommand{\ecor}{\end{corr}}
\newcommand{\edeff}{\end{deff}}
\newcommand{\ii}{i}
\newcommand{\eprop}{\end{proposition}}

\newcommand{\supp}{\text{supp }}
\renewcommand{\Pi}{\varPi}

\renewcommand{\epsilon}{\varepsilon}

\newcommand{\Adel}{{{\cal A}}_\delta}

\begin{document}

\title[Concerning Nikodym-type sets in curved spaces]
{Concerning Nikodym-type sets in 3-dimensional curved spaces}
\thanks{
MR Classification Numbers: 42, 58.
\\
The author was supported in part by the NSF}
\author{Christopher D. Sogge}
\address{Department of Mathematics, The Johns Hopkins University, Baltimore,
Maryland 21218}
\email{sogge\@jhu.edu}

\begin{abstract}
We investigate maximal functions involving averages over geodesics in three-dimensional Riemannian manifolds.  We first show that one can easily extend the Euclidean results of Bourgain and Wolff if one assumes {\it constant curvature}.  If this assumption may not hold.  Nonetheless, we formulate a generic geometric condition which allows favorable estimates.  Curiously, this condition ensures that one is in some sense as far as possible from the constant curvature case.  Assuming it one can prove dimensional estimates for Nikodym-type sets which are essentially optimal.  Optimal estimates for the related maximal functions are still open though.
\end{abstract}

\maketitle

\newsection{Introduction}

In this paper we shall give some natural partial extensions to the curved space setting of results of Bourgain \cite{B1} and Wolff \cite{W} concerning lower bounds for the dimension of compliments of Nikodym sets in Euclidean space.

Recall that a classical Nikodym set is a subset of $[-1,1]\times [-1,1]$ of Lebesgue measure one which has the property that for each $x\in F$ there is a line $\gamma_x$ so that $\gamma_x \cap F=  \{x\}$.  Because of this, the relative compliment, $\Omega = [-1,1]\times [-1,1]\backslash F$, must be a set of measure zero with the property that if $\Omega^*_\alpha$ is the set of points $x$ for which there is a line segment $\gamma_x^\alpha$ through $x$ with $|\Omega \cap \gamma^\alpha_x| =\alpha$, then $|\Omega^*_\alpha|>0$ for every $0<\alpha<1$. Here $|\Omega\cap \gamma^\alpha_x|$ denotes one-dimensional Lebesgue measure.

Results of C\'ordoba \cite{C} imply that such a set must have full Hausdorff dimension.  For analogous sets in ${\Bbb R}^3$ it is conjectured that the same should be true.  By taking projections, the results of \cite{C} immediately imply that such sets must have dimension at least 2.

This result was first improved by Bourgain \cite{B1}.  His results say that if $\Omega\subset {\Bbb R}^3$ and if $\Omega^*_\alpha=\{x: \, |\Omega \cap \gamma^\alpha_x|=\alpha \}$ then $\text{dim }\Omega \ge 7/3$ if $|\Omega^*_\alpha|>0$.  Here, as before, $\gamma^\alpha_x$ denotes a line segment through $x$ of length $\alpha>0$.  This lower bound was later improved by Wolff \cite{W} to $\text{dim }\Omega\ge 5/2$ if $|\Omega^*_\alpha|>0$.  In both works the lower bounds on the dimension were obtained for somewhat more general sets.  Lower bounds for analogous sets in higher dimensions were also obtained in \cite{B1} and \cite{W}.  The strongest to date are those of Wolff \cite{W} showing that if $\Omega \subset {\Bbb R}^n$ and $|\Omega^*_\alpha|>0$ for some $\alpha>0$ then $\text{dim }\Omega \ge (n+2)/2$.

Let us now consider extensions of this result to the curved $3$-dimensional setting.  To this end, we shall let $(M^3,g)$ denote a (paracompact) $3$-dimensional Riemannian manifold with metric $ds^2=\sum g_{jk}(x) dx_j dx_k$.  Given $x\in M^3$, we let $\{\gamma^\alpha_x\}$ denote the set of all geodesics containing $x$ with arclength $\alpha$, that is, $|\gamma_x^\alpha|=\alpha$.  Abusing the classical terminology somewhat, we now define Nikodym-type subsets of $M^3$.

\bdeff\label{def1.1}
If $\Omega\subset M^3$, $\alpha>0$ and $0<\lambda<1$, let 
\begin{equation}\label{1.1}
\Omega_{\alpha,\lambda}^* =\{x\in M^3: \, \exists\,  \gamma^\alpha_x \, \, \text{with } |\gamma^\alpha_x\cap\Omega|\ge \lambda |\gamma^\alpha_x|\}.
\end{equation}
We then say that $\Omega$ is a {\it Nikodym-type} set if, for 
a finite $\alpha$ smaller than half the injectivity radius of $M^3$
and all $\lambda$ sufficiently close to $1$, $\Omega^*_{\alpha,\lambda}$ has positive measure.
\edeff

A couple of remarks are in order.  First, if $M^3$ is Euclidean space ${\Bbb R}^3$, then these sets are slightly more general than the ones mentioned before.  Nonetheless, the lower bounds mentioned before of Bourgain and Wolff hold for the Euclidean case if one just assumes that $|\Omega^*_{\alpha,\lambda}|>0$ for some $\alpha>0$ and $0<\lambda<1$.

We shall see that Wolff's lower bound $\text{dim }\Omega\ge 5/2$ holds if $(M^3,g)$ has {\it constant}  curvature.  The proof merely involves a straightforward adaptation of Wolff's argument using Fermi normal coordinates.  The only minor difference in the main part of our argument versus that in \cite{W} is that we rely on $L^2$-bounds for a weighted auxiliary maximal function.  This fortunately allows us to avoid the inductive argument in \cite{W} which relied on a simple scaling argument which seems difficult to generalize to the non-Euclidean setting. 

The arguments involved rely on the fact that if $M^3$ has constant curvature and if Fermi (local) coordinates are chosen about a geodesic segment $\gamma$ then every resulting local ``Fermi two-plane" is totally geodesic.  This fact and the argument that exploits it of course are not stable under perturbations. 

Based on this principle it was shown by Minicozzi and the author \cite{MS} that for general Riemannian manifolds $M^3$ the ``easy" lower bound $\text{dim }\Omega \ge 2$ for Nikodym-type sets is in general sharp even if $\text{dim }\Omega$ refers to the Minkowski dimension.  Indeed if, for any $\varepsilon>0$, we consider ${\Bbb R}^3$ with the metric 
$dx^2+\varepsilon a(x_1)dx_2dx_3$ where $a(s)=e^{1/s}, s<0$ and $a(s)=0, s\ge 0$, then the subset $\Omega =\{x: \, x_3=0, \, 0<x_1, \, |x_2|\le 1\}$ of the two-plane where $x_3=0$ is a Nikodym-type set. This is because there is a neighborhood ${\cal N}$ of $\{(x_1,0,0): x_1<0\}$ so that if $x\in {\cal N}$, then there is a geodesic $\gamma_x\ni x$ which lies in the two plane where $x_3=0$ if $x_1>0$ and intersects $\Omega$ in a set of positive measure. In this example, all sectional curvatures vanish when $x_3\ge0$.   Similar considerations show that if $(M^3,g_0)$ has constant non-zero curvature then one can find an arbitrarily small perturbation $(M^3,g)$ so that the resulting Riemannian manifold $(M^3,g)$ has Nikodym-type sets with Minkowski dimension two.

Keeping this in mind, if one wishes to improve on the lower bounds $\text{dim }\Omega \ge 2$ in \cite{MS}, it is necessary to place a condition on $(M^3,g)$ to ensure that one is quantitatively ``far" from the constant curvature case.  In three dimensions such a condition is easy to formulate using the Ricci tensor.  Under this curvature hypothesis on $(M^3,g)$ we shall show that $\Omega$ must have Minkowski dimension at least $7/3$ if $\Omega\subset M^3$ is a Nikodym-type set, recovering an analog of the lower bound of Bourgain \cite{B1} for this case.  Obtaining Wolff's \cite{W} lower bound $\text{dim }\Omega\ge5/2$ seems much more delicate here.  In particular, at the end of this paper, we shall see that the $L^{5/2}\to L^{10/3}$ maximal estimates on which the Euclidean lower bounds are based cannot hold in ``variably curved" manifolds (see Definition 3.1).  We shall see that in closely related cases involving the more local ``chaotic curvature condition'' (see Definition 3.1), the best one could hope for would be the analog of the Euclidean $L^{7/3}\to L^{7/3}$ estimate of Bourgain \cite{B1}.  However, for the more typical ``variably curved'' case a substitute for Wolff's result may hold if one replaces $L^{10/3}$ by $L^{5/2}$.  As we shall see at the end, this would be sharp.
Basically, the bounds for the maximal functions must be worse than those in the Euclidean case due too small-scale ``focusing" of geodesics while this might not be an obstacle for the problems involving lower bounds for the dimensions of Nikodym-type sets.

The paper is organized as follows.  In the next section we shall prove the results mentioned earlier for spaces of constant curvature.  Here we shall also see that Wolff's bounds for Euclidean Nikodym maximal functions extend easily in this case if one considers tubes whose length $\alpha$ is smaller than half the injectivity radius. Such estimates of course easily give the lower bounds for the dimension of Nikodym-type sets, and, as in the Euclidean case, we only need to assume that $|\Omega^*_{\alpha,\lambda}|>0$ for some $\lambda>0$ to conclude that $\text{dim }\Omega\ge 5/2$.  In Section 5 we shall see that for ``variably curved" manifolds all of the arguments for the constant curvature case break down due to the fact that auxiliary maximal functions involving averaging over small tubes about geodesics intersecting a common geodesic have unfavorable $L^2\to L^2$ bounds.  Despite this, in Section 3 we shall be able to obtain our lower bound $\text{dim }\Omega \ge 7/3$ for Nikodym-type sets using considerably weaker auxiliary $L^2$-estimates.  At present, we unfortunately do not have what seem to be natural related $L^{7/3}$-estimates for the Nikodym maximal functions in this context, because, in part, of the difficulty in dealing with small scales where the geometry becomes Euclidean.

Let us say a few words about the auxiliary $L^2$ estimates which we shall employ since this is crux of our analysis.  As mentioned before, in all cases, the maximal functions involved in them just involve averages over thin tubes whose centers intersect a fixed common geodesic $\gamma_0$.  In the constant curvature case, it is fairly easy to prove ``sharp'' estimates for this operator since we can reduce matters to a simple two-dimensional estimate if we just use the fact that in (local) Fermi normal coordinates about $\gamma_0$ every two plane containing $\gamma_0$ is totally geodesic, and therefore the center of any tube in our average must be contained in one of these two planes. In the variably curved case this argument of course completely breaks down.  Here our $L^2$ estimates are based on properties of the Fourier integrals underlying the averages.  We shall exploit the fact that under our geometric assumptions they have canonical relations that, off possibly a small exceptional set, have projections with at most folding singularities in the sense of Melrose and Taylor \cite{MT}.  The reason for this is that in the ``variably curved'' case, unlike the constant curvature case, the set of geodesics intersecting a common geodesic $\gamma_0$ is fairly randomly distributed.  We should point out that in the constant curvature case the underlying Fourier integral operators are much more degenerate; however, this is more than compensated by their concentration properties mentioned before.  To exploit the salient features of the Fourier integrals governing the auxiliary averages in the ``variably curved'' case we use a theorem of Melrose and Taylor \cite{MT} which says that Fourier integral operators with folding canonical relations are bounded on $L^2$ with a loss of $1/6$ of a derivative.  If we apply their theorem we lose $1/6$ of a power of $\delta$ in our auxiliary maximal function bounds and this accounts for our lower bound of $7/3=5/2-1/6$ versus $5/2$ for the dimension of Nikodym-type sets.  The fact that we have to avoid an exceptional set where the Fourier integrals may be more degenerate accounts for our assumption that $|\Omega^*_{\alpha,\lambda}|>0$ for all $O<\lambda<1$ in Definition 1.1.  On the other hand, it turns out that our Fourier integral operators have the property that at least one of the projections of the canonical relation has at most folding singularities.  In Section 4 we shall exploit this fact and appeal to a theorem of Greenleaf and Seeger \cite{GS} which says that such operators are always bounded on $L^2$ with a loss of $1/4$ of a derivative.  By doing so we shall be able to prove slightly less favorable lower bounds for a wider class of sets in ``variably curved'' manifolds.  Specifically, if $\alpha$ is small and if $|\Omega^*_{\alpha,\lambda}|>0$ for some $0<\lambda<1$, then we shall be able to show that $\text{dim }\Omega\ge 9/4$.

Throughout this paper $C$ and $c$ will denote positive finite constants which may change at each occurrence. Also, to avoid burdensome notation, we shall be inconsistent by using in different places subscripts and superscripts to denote local coordinates; however, the meaning should be clear in the given context.

It is a pleasure to thank my colleagues W. Minicozzi, J. Spruck and S. Zelditch for many helpful patient discussions and suggestions.  I would also like to thank A. Seeger for a helpful discussion regarding general x-ray transforms.  The author also benefited from a course taught by M. Christ on Wolff's paper \cite{W}.

\newsection{Spaces of constant curvature}

Let $(M^3,g)$ be a Riemannian manifold, and, as before, let $\gamma^\alpha_x$ denote all geodesics passing through $x\in M^3$ of length $|\gamma^\alpha_x|=\alpha$.  We assume that $\alpha>0$ is finite and smaller than half of the injectivity radius of $M^3$.  Using the metric, we then let 
$$T^{\alpha,\delta}_x=\{y\in M^3: \, \text{dist }(y, \gamma^\alpha_x)\le \delta\},$$
be a tubular neighborhood of width $\delta$ around $\gamma^\alpha_x$.  We shall also at times slightly change the notation, denoting the tubes for instance by $T^{\alpha,\delta}_{\gamma^\alpha_x}$. Given a function $f$, we can now define the Nikodym maximal functions
\begin{equation}\label{nikmax}
f^*_\delta(x)=\sup |T^{\alpha,\delta}_x|^{-1}\int_{T_x^{\alpha,\delta}} |f(y)|\, dy.
\end{equation}

If $M^3$ is flat Euclidean space ${\Bbb R}^3$, Wolff \cite{W} obtained the following estimates for these operators when $1\le p\le 5/2$:
\begin{equation}\label{wolffbound}
\|f^*_\delta\|_{L^{2p/(p-1)}({\Bbb R}^3)}\le C_\varepsilon \delta^{1-3/p-\varepsilon}\|f\|_{L^p({\Bbb R}^3)}, \, \, 1\le p\le 5/2, \, \, \varepsilon>0.
\end{equation}
Except for the $\varepsilon$, these bounds are easily seen to be best possible.

The main result of this section is that bounds like this hold in the constant curvature case.

\bth\label{thm2.1}
Assume that $(M^3,g)$ has constant curvature.  Then for $f$ supported in a compact subset $K$ of a coordinate patch and all $\varepsilon>0$
\begin{equation}\label{constant}
\|f^*_\delta\|_{L^{2p/(p-1)}(M^3)}\le C_\varepsilon \delta^{1-3/p-\varepsilon}\|f\|_{L^p(M^3)}, \, \, 1\le p\le 5/2, \, \, \supp f\subset K.
\end{equation}
\eth

Since $|T^{\alpha,\delta}_x|\approx \delta^2$, the $L^1\to L^\infty$ bounds are trivial.  The other estimates follow via interpolation from a restricted weak-type inequality corresponding to $p=5/2$.  This says that if $f=\chi_E$ is the characteristic function of a set $E\subset K$, then, for all $\varepsilon>0$,
\begin{equation}\label{weak}
|\{x: \, f^*_\delta(x)\ge \lambda\}|\le C_\varepsilon \bigl(\, \delta^{-1/2-\varepsilon}\lambda^{-5/2}|E|\, \bigr)^{4/3}, \, \, f=\chi_E.
\end{equation}
We should point out that the proof below also yields the sharper result where in \eqref{weak} $\delta^{-1/2-\varepsilon}$ is replaced by $\delta^{-1/2}$ times a sufficiently high power of $\log\delta$ if $1\le p<5/2$.  Such a refinement, though,  would complicate the bookkeeping a bit, and  is not important for the applications we have in mind.

If we use this estimate with $\lambda$ close to one (as in Definition 1.1) we can show that the lower bounds in \cite{W} on the Minkowski dimension\footnote{
Similar arguments using \eqref{weak} for $\lambda\to 0_+$ imply stronger results involving the same lower bounds for the Hausdorff dimensions (see \cite{B1}, \cite{W}).} of Nikodym-type sets in Euclidean space carry over to the present setting.

\bcor\label{cor2.1}
If $\Omega$ is a Nikodym-type set in a manifold $M^3$ of constant curvature, then
\begin{equation}\label{2.5}
\text{dim }\Omega \ge 5/2.
\end{equation}
Moreover, if $|\cup_{\lambda>0}\Omega^*_{\alpha,\lambda}|>0$ then the same conclusion holds if $\Omega^*_{\alpha,\lambda}$ is as in \eqref{1.1}.
\ecor

\begin{pf}
There is no loss of generality in assuming that $\Omega\subset K$, where, as above, $K$ is a compact subset of a coordinate patch in $M^3$.  Let $\Omega_\delta=\{x: \, \text{dist }(x,\Omega)\le \delta \, \}$ be a $\delta$-neighborhood of $\Omega$.  We then must show that if $\varepsilon>0$ is fixed there is a uniform constant $c_\varepsilon>0$ so that
\begin{equation}\label{minconst}
|\Omega_\delta|\ge c_\varepsilon \delta^{3-5/2+\varepsilon}, \, \, 0<\delta<1,
\end{equation}
assuming that $|\Omega^*_{\alpha,\lambda}|>0$, where $\Omega^*_{\alpha,\lambda}$ is as in Definition 1.1.

The proof of \eqref{minconst} is easy.  If we take $f=\chi_{\Omega_\delta}$ in \eqref{weak}, then $f^*_{C\delta}(x)\ge \lambda/C$ when $x\in \Omega^*_{\alpha,\lambda}$ if $C$ is a large fixed constant.  Hence \eqref{weak} yields
$$|\Omega_\delta|\ge C_\varepsilon'\delta^{1/2+\varepsilon}\, |\Omega_{\alpha,\lambda}^*|^{3/4},$$
leading to \eqref{minconst} since we are assuming that $|\Omega^*_{\alpha,\lambda}|>0$ for some $\lambda>0$.
\end{pf}

Let us turn to the proof of \eqref{weak}.  The key ingredient is an $L^2$-estimate for an auxiliary maximal operator.  As we shall see this estimate is what breaks down if one does not assume constant curvature.

The maximal operator involves averages over small neighborhoods of geodesics $\gamma_x\ni x$ which intersect a  {\it fixed} ``common'' geodesic $\gamma_0$.  Here and throughout the rest of the paper, we shall assume that all the geodesics involved have length $\alpha$, where $\alpha$ is assumed to be fixed and  to be no more than the minimum of  $1$ and half of the injectivity radius of $M^3$.  We shall also drop the various dependencies on $\alpha$ from the notation.

We could use a coordinate-free definition of our auxiliary maximal operator; however, for the proof of its bounds, and for the application, it is convenient to work in a special (local) coordinate system about $\gamma_0$ called Fermi normal coordinates.  Recall that these amount to a generalization of geodesic normal coordinates where a point is replaced by a geodesic.  (See, e.g., \cite{HE} Section 4.1, \cite{H2} and \cite{tubesbook} for further discussion.)

One obtains Fermi normal coordinates in the following manner.  First, one fixes a point $x_0\in \gamma_0$ and then chooses an orthonormal basis $E_1,E_2,E_3 \in T_{x_0}M^3$ with $E_1$ being a unit tangent vector of $\gamma_0$ at $x_0$.  Using parallel transport, one propagates this basis to every point of $\gamma_0$.  If $\gamma_0$ is the arclength parameterization of $\gamma_0$ with $\gamma_0(0)=x_0$, $\dot \gamma_0(0)=E_1$, then the resulting vectors $E_j(s)$ are orthonormal in $T_{\gamma_0(s)}M^3$ and $\dot \gamma_0(s)=E_1(s)$.  One then assigns Fermi coordinates $(x_1,x_2,x_3)$ to a point $x$ if $x$ lies a distance $|(x_2,x_3)|$ from $\gamma_0$ on the geodesic through $\gamma_0(x_1)$ whose unit tangent vector at this point is $(x_2E_2(x_1)+x_3E_3(x_1))/|(x_2,x_3)|$.

These coordinates are well defined near $\gamma_0$.  Note further that in these coordinates the metric must satisfy
\begin{equation}\label{fermi1}
\sum_{k=2}^3g_{jk}(x)x_k=
\cases x_j, \quad  \text{if } \, 2\le j\le 3
\\
0, \quad \text{if } \, j=1,
\endcases
\end{equation}
as well as
\begin{equation}\label{fermi2}
(\partial/\partial x_2)^m (\partial/\partial x_3)^n
(g_{jk}(x)-\delta_{jk})=0 \, \text{ if } m+n<2 \, \, \text{ and }\, x_2=x_3=0.
\end{equation}
The first condition means precisely that the rays $t\to (x_1,tx_2,tx_3)$ are geodesics orthogonal to $\gamma_0=\{(s,0,0)\}$.  The second follows from the first and the fact that for every $j$ $\partial/\partial x_j$ is parallel along $\gamma_0$.  Also, note that these Fermi normal coordinates are unique up to rotations preserving the $x_1$-axis.

We can now define our auxiliary maximal function using these coordinates.  If $x'=(x_2,x_3)$, we set
\begin{equation}\label{auxdef}
(\Adel f)(x') =
\sup_{\{\gamma_{x'}\ni (0,x'): \, \gamma_{x'}\cap\gamma_0\ne \emptyset \}}
|T^\delta_{\gamma_{x'}}|^{-1}
\int_{T^\delta_{\gamma_{x'}}} |f(y)|\, w_{\gamma_{x'}}(y)\, dy,
\end{equation}
where the damping factor is
\begin{equation}\label{damping}
w_{\gamma_{x'}}(y)=\bigl(\, \text{dist }(y,\gamma_{x'}\cap\gamma_0)\, \bigr)^{1/2} \, .
\end{equation}
Thus $(\Adel f)(x')$ should be thought of as a maximal function that takes its values on a hypersurface which is transverse to $\gamma_0$, and it just involves averages over $\delta$-neighborhoods of geodesics intersecting the common geodesic $\gamma_0$.

The estimates we require of $\Adel$, which are essentially of the best possible nature, are the following

\ble\label{auxest1}  Assume that $M^3$ has constant  curvature.  Then for $f$ supported in a compact subset $K$ of a coordinate patch \begin{equation}\label{good}
\|\Adel f\|_{L^2}\le C (\log 1/\delta)^{3/2} \|f\|_{L^2}, \, \, 0<\delta<1/2.
\end{equation}
\ele

In the Euclidean case, this is just a minor modification of Lemma 2.1 in Wolff \cite{W}.  In the present context, \eqref{good} is a simple consequence of a variable coefficient version of a maximal theorem of C\'ordoba \cite{C} (see also \cite{MSS}) involving averages of functions of two variables.  We postpone the straightforward argument until the end of this section. 

Using this lemma we can prove \eqref{weak} using multiplicity arguments as in Wolff \cite{W}.  First, though, as in \cite{B1} or \cite{W}, it is convenient to state a discrete form of the problem.

To do this, we first use the induced metric on the unit tangent bundle to define the 
$TM^3$-distance between two geodesics $\gamma_j$, $j=1,2$, of length $\alpha$.  Specifically, we put
\begin{equation}\label{tdist}
\text{dist}_{TM^3}(\gamma_1,\gamma_2) = \min_{x_j\in \gamma_j, \tau_j=\dot \gamma_j|_{\gamma_j=x_j}} \text{dist } ((x_1,\tau_1), (x_2,\tau_2)).
\end{equation}
Here $\dot\gamma_j|_{\gamma_j=x_j}$ denotes a unit tangent vector at $x_j$.

If we fix a geodesic $\gamma_0$, we then consider the family of all geodesics where
\begin{equation}\label{small}
\text{dist}_{TM^3}(\gamma,\gamma_0)\le c,
\end{equation}
with $c>0$ being a small fixed constant.  Working in the Fermi normal coordinates about $\gamma_0$, we further assume that $\gamma$ intersects the hyperplane $\{(0,x'): \, x'=(x_2,x_3)\}$; we do not, though, of course assume now that $\gamma$ intersects $\gamma_0$.  Let us call the resulting family of geodesics
\begin{equation}\label{close}
{\cal F}= \{\gamma_{x'}: \, (0,x')\in \gamma_{x'} \, \text{and } \, \text{dist}_{TM^3}(\gamma_{x'},\gamma_0)\le c\}.
\end{equation}
(Note that over every $x'$ there is a two-parameter family of geodesics $\gamma_{x'}\in {\cal F}$.)  We then consider a $C_0\delta$-separated collection of points
\begin{equation}\label{basepoint}
\{x'_j: \, 1\le j\le M\},
\end{equation}
where $C_0$ is a fixed constant, 
and choose, for each $j$, a geodesic $\gamma_{x_j}\in {\cal F}$.  If
we assume further that for some fixed $\lambda\in [\delta,1]$ 
\begin{equation}\label{assumption}
|E\cap T^\delta_j|\ge \lambda |T^\delta_j|, \quad 1\le j\le M
\end{equation}
where 
$$T^\delta_j=\{y: \, \text{dist }(y,\gamma_{x_j})\le \delta\}$$
is the $\delta$-tube about $\gamma_j$, then \eqref{weak} would follow from the uniform bounds
\begin{equation}\label{discrete}
M\delta^2 \le C_\varepsilon \bigl(\delta^{-1/2-\varepsilon}\lambda^{-5/2}|E|\bigr)^{4/3}.
\end{equation}
Indeed, this inequality is equivalent to the slightly stronger version of \eqref{weak} where we replace the left side by $|\{x': \, f^*_\delta(0,x')\ge \lambda\}|,$ $f=\chi_E$ and replace the maximal operator by one involving averaging over $\delta$-tubes with centers satisfying \eqref{close}.  Note also, for later use that since the basepoints $\{x'_j\}$ of the tubes are $\delta$-separated, we must have 
\begin{equation}\label{angle}
\text{angle }(T^\delta_{i}, T^\delta_j) \ge c\delta\quad
\text{if } \, T^\delta_i\cap T^\delta_j\ne \emptyset,
\end{equation}
for some uniform $c>0$, if $\text{angle }(T^\delta_i, T^\delta_j)
=\min \text{dist }(\dot\gamma_i, \dot \gamma_j),$ where the minimum is taken over points in a $\delta$-ball of $T^i_\delta\cap T^\delta_j$.   To simplify the notation, we can assume that $c$ is a large fixed constant if we fix $C_0$ above large enough.

To proceed, we shall use a slight variation of the multiplicity argument in \cite{W}.  Our modification will allow us to avoid the induction argument on the eccentricity of the tubes in \cite{W}, which is fortunate since ``scaling'' arguments are much more complicated in the non-Euclidean setting.

Let us be more specific.  First, if $1\le j\le M$ and $x\in T^\delta_j$ are fixed, let
\begin{equation}\label{anglemult}
{\cal I}_\theta(x,j)=\{i: \, x\in T^\delta_i \, \, \, \text{and } \, \text{angle } (T^\delta_i, T^\delta_j)\in [\theta/2, \theta] \}
\end{equation}
index the tubes $T^\delta_i$  intersecting $T^\delta_j$ at $x$ with angle $\approx \theta$.  Next, let
\begin{multline}\label{shellmult}
{\cal I}_\mu(x,j)=
\\
\{i:  x\in T^\delta_i \, \,  \text{and } \,  
|T^\delta_i \cap \{y\in E: \, \text{dist }(y,\gamma_j)\in [\mu/2, \mu]\}|\ge 
(2\log_2 1/\delta)^{-1}\lambda |T^\delta_i|\}
\end{multline}
index the tubes $T^\delta_i$ which intersect $T^\delta_j$ at $x$ and have the property that a non-trivial portion of $E$ belongs to the part of $T^\delta_i$ which is a distance $\approx \mu$ away from $T^\delta_j$.

If we combine these two conditions, videlicet, 
$${\cal I}_{\theta\mu}(x,j)={\cal I}_\theta(x,j)\cap {\cal I}_\mu(x,j),$$
then we have the following

\ble\label{multlemma}
There 
are
$N\in {\Bbb N}$ and $\theta,\mu\in [\delta,1]$
so that there are at least $M/2$ values of $j$ for which 
\begin{equation}\label{case1}
|\{x\in T^\delta_j\cap E: \, \text{card }(\{i: x\in T^\delta_i \})\le N\}| \ge (\lambda/2) |T^\delta_j|,
\end{equation}
and, moreover, 
\begin{equation}\label{case2}
|\{x\in T^\delta_j\cap E: \, \text{card }{\cal I}_{\theta\mu}(x,j)\ge N/(2\log_2 1/\delta)^2\}| \ge (4\log_2 1/\delta)^{-2}\lambda |T^\delta_j|
\end{equation}
for at least $M/(2\log_21/\delta)^2$ values of $j$.
\ele

\begin{pf}  Choose the smallest $N\in {\Bbb N}$ so that \eqref{case1} holds.  Then there must be $M/2$ values of $j$ for which
\begin{equation}\label{blah}
|\{x\in T^\delta_j\cap E: \, \, \text{card }(\{i: x\in T^\delta_i\})\ge N\}| \ge (\lambda/2)|T^\delta_j|.
\end{equation}
For any such fixed $j$ and $x\in T^\delta_j\cap E$ with $\text{card }(\{i: x\in T^\delta_i\})\ge N$ we can find $1\le m,n\le \log_2\delta$ so that
$${\cal I}_{2^m\delta,2^n\delta}(x,j)\ge N/(2\log_2 1/\delta).$$
To verify this one uses \eqref{assumption}, \eqref{angle} and our assumptions that, in these two inequalities, $\lambda/\delta$ and $c$ are large.  Finally, since there are $M/2$ values of $j$ satisfying \eqref{blah}, if we use the pigeonhole principle one more time, we conclude that we can choose fixed $\theta=2^m\delta$ and $\mu=2^n\delta$ so that \eqref{case2} holds for at least $M/(2\log_2 1/\delta)^2$ values of $j$, which finishes the proof.
\end{pf}

Given this splitting, the main step in the proof of \eqref{discrete} is to obtain the following

\bprop\label{mainprop}  Let $N$ be as in the preceding lemma.  Then
\begin{equation}\label{trivial}
|E|\ge \lambda M\delta^2/CN.
\end{equation}
Moreover, if $T^\delta_j$ is a tube for which \eqref{case2} holds and if $\varepsilon>0$ is fixed there is a uniform constant $C_\varepsilon$ so that, for small $\delta>0$, given $a\in M^3$
\begin{equation}\label{nontrivial}
|(E\backslash B(a,\delta^\varepsilon\lambda) )\cap T^\mu_j|\ge C_\varepsilon \lambda^3\mu \delta^{1+\varepsilon}N,
\end{equation}
if $B(a,r)=\{y\in M^3: \text{dist }(y,a)\le r\}$.
\eprop

Before turning to the proof, let us see how these two estimates lead to \eqref{discrete}.  

The first step is to realize that we can use \eqref{nontrivial} and Bourgain's bush argument in \cite{B1} to obtain another lower bound of $|E|$ involving $M$.  Specifically, we claim that
\begin{equation}\label{punch}
|E|\ge c\lambda^4 N\delta^{1+2\varepsilon}\sqrt{M\delta^2}.
\end{equation}
Clearly this inequality and \eqref{trivial} imply \eqref{discrete} if one takes geometric means.

To verify \eqref{punch} we require the following

\ble\label{bushlemma}
Suppose that $T^\mu_j$, $1\le j\le M_0$ are tubes of thickness and length $\alpha$, where $\alpha$ is as above.  Assume also that for some $0<\theta<1$
\begin{equation}\label{bass1}
\text{angle }(T^\mu_j,T^\mu_k)\ge C_1\theta \quad \text{if } \, T^\mu_j\cap T^\mu_k\ne \emptyset.
\end{equation}
Assume also that for every $a\in M^3$
\begin{equation}\label{bass2}
|T^\mu_j\cap (E\backslash B(a,\mu/\theta))|\ge \rho |T^\mu_j|, \quad
1\le j\le M_0.
\end{equation}
Then if $C_1$ is large enough, there is a fixed $c>0$ so that
\begin{equation}\label{bconcl}
|E|\ge c\rho\mu^2\sqrt{M_0}.
\end{equation}
\ele

We shall postpone the proof. For now, let us see why it along with \eqref{nontrivial} leads to \eqref{punch}.  To do this, if ${\cal J}$ denotes the $M/(2\log_2 1/\delta)^2$ values of $j$ for which \eqref{nontrivial} holds, and if $\{x'_j\}$ are the corresponding basepoints as in \eqref{basepoint}, let us choose a subcollection $\{j_k\}_{k=1}^{M_0}\subset{\cal J}$ so that the resulting points $x'_{j_k}$ are $C_0(\mu/\delta^\varepsilon\lambda)$-separated.   Using the inclusion relation, if we choose a maximal such subcollection we must have
$$M_0\ge c M(\delta/(\mu/\delta^\varepsilon\lambda))^2=cM(\lambda\delta^{1+\varepsilon}/\mu)^2,$$
where $c$ depends on $C_0$.  If the latter constant is large, the associated tubes $T^\mu_{j_k}$ verify \eqref{bass1} with $\theta= \mu/\delta^\varepsilon\lambda$.  Since $\mu/\theta=\delta^\varepsilon\lambda$, by \eqref{nontrivial}, we must have \eqref{bass2} with 
$$\rho\approx \lambda^3N\delta^{1+\varepsilon}\mu/\mu^2=\lambda^3N\delta^{1+\varepsilon}/\mu.$$
Thus, since $\rho\mu^2\sqrt{M_0}$ must be larger than a multiple of
$$\lambda^3N\delta^{1+\varepsilon}\mu^{-1}\cdot \mu^2\cdot \sqrt{M(\lambda\delta^{1+\varepsilon}/\mu)^2}=\lambda^4N\delta^{1+2\varepsilon}
\sqrt{M\delta^2},
$$
\eqref{punch} follows from \eqref{bconcl} and \eqref{nontrivial} as claimed.

We now turn to the proof of  Proposition \ref{mainprop}.

\noindent{\it Proof of \eqref{trivial}.}  If we let $E_0=\{x\in E: \, \sum_{k=1}^M \chi_{T^\delta_k(x)}\le N\}$, then, by the first part of Lemma \ref{multlemma}, $|T^\delta_j\cap E_0|\ge \lambda|T^\delta_j|/2$ for at least $M/2$ values $j=j_k$.  Thus, since $|T^\delta_{j_k}|\approx \delta^2$,
$$|E|\ge
|\cup(E_0\cap T^\delta_{j_k})|
\ge N^{-1}
\sum_{k=1}^{M/2}|E_0\cap T^\delta_{j_k}|\ge c\lambda M\delta^2/N,$$
as desired.

\medskip

\noindent{\it Proof of \eqref{nontrivial}.}  Fix $j$ as in \eqref{case2}.  Then if $i\in {\cal I}_{\theta\mu}$ recall that $T^\delta_i$ intersects $T^\delta_j$ at angle $\approx \theta$ and that
\begin{equation}\label{1}
|T^\delta_i \cap \{y\in E: \, \text{dist }(y,\gamma_j)\in [\mu/2, \mu]\}|\ge 
(2\log_2 1/\delta)^{-1}\lambda |T^\delta_i|.
\end{equation}
Since $|T^\delta_i \cap B(a,\delta^\varepsilon\lambda)|\le C\delta^\varepsilon \lambda |T^\delta_i|$, $a\in M^3$, if we replace $E$ in \eqref{1} by $E\backslash B(a,\delta^\varepsilon\lambda)$ we have the same sort of lower bound if the $2$ in the right is replaced by $4$ when $\delta$ is small.  Hence, if we replace $\lambda$ by $\lambda/2$, we conclude that \eqref{nontrivial} would follow if we could show that if \eqref{case2} holds then
\begin{equation}\label{newmain}
|E\cap T^\mu_j|
\ge C_\varepsilon \lambda^3\mu\delta^{1+\varepsilon}N.
\end{equation}

To prove this we shall use Lemma \ref{auxest1}.  If we let $\gamma_0$ in there be the center $\gamma_j$ of $T^\delta_j$ and work in Fermi normal coordinates about this geodesic, we need a localized discrete form of \eqref{good}.  Specifically, let $x'_i$, $1\le i\le M_0$ be a $\delta$-separated collection of basepoints and assume that for every $i$ that there is a tube $T^\delta_i$ containing $(0,x'_i)$ which intersects the fixed common tube $T^\delta_j$.  Assume further that
\begin{equation}\label{a}
\text{angle }(T^\delta_j,T^\delta_i)\approx \theta, \quad 1\le i\le M_0
\end{equation}
and that
\begin{equation}\label{0}
|T^\delta_i \cap \{y\in E: \, \text{dist }(y,\gamma_j)\in [\mu/2,\mu]\}|\ge \rho, \quad 1\le i\le M_0.
\end{equation}
Note that the preceding inequality yields lower bounds for ${\cal A}_{C\delta}\chi_E(x'_i)$, for sufficiently large $C$, since $T^\delta_i$ intersects $T^\delta_j$.  With this in mind, we claim that \eqref{good} along with \eqref{a} and \eqref{0} yield
\begin{equation}\label{2}
|E|\ge C(\mu/\theta)M_0\delta^2\rho^2/(\log 1/\delta)^3.
\end{equation}

To verify this, we first note that \eqref{good} of course implies the corresponding weak-type bounds
\begin{equation}\label{3}
|\{x': \, \Adel \chi_E(x')>\rho\}|\le C\rho^{-2}(\log 1/\delta)^3|E|.
\end{equation}
Next, since \eqref{a} implies the damping factors $w_{\gamma_{x'_i}}(y)$ in the definition \eqref{auxdef} of $\Adel$ are $\approx (\mu/\theta)^{1/2}$ on the set in the left side of \eqref{0}, we conclude that we must have lower bounds of the form
\begin{equation}\label{pig}
({\cal A}_{C\delta} f)(x'_i)\ge c(\mu/\theta)^{1/2}\rho,\quad 1\le i\le M_0
\end{equation}
for uniform $C$ and $c>0$ if \eqref{a} and \eqref{0} hold.  
We then obtain the discrete inequality \eqref{2} from \eqref{pig} in the same way that \eqref{discrete} follows from its corresponding weak-type inequality. 

For the next step, we claim that there must be at least 
\begin{equation}\label{4}
M_0\ge cN\lambda\theta/\delta(\log 1/\delta)^4
\end{equation}
tubes $T^\delta_{\ii_k}$ which intersect $T^\delta_j$ and satisfy \eqref{1} and \eqref{a}.  If we then take $\rho=\lambda/(2\log 1/\delta)$, we obtain \eqref{newmain} from \eqref{2}.

The proof of \eqref{4} is straightforward.  By \eqref{case2} we must have
$$\sum_{k=1}^{M_0}\chi_{T^\delta_{\ii_k}}(x)\ge N/(2\log 1/\delta)^2$$
when $x$ belongs to a subset of $T^\delta_j$ of measure $\lambda|T^\delta_j|/(4\log 1/\delta)^2$.  Note further that since $T^\delta_j$ and such a $T^\delta_{\ii_k}$ intersect at angle $\approx \theta$, we have
$|T^\delta_j\cap T^\delta_{\ii_k}|\ge c\delta^3/\theta$ (see Lemma \ref{inj} below).  Putting all of this together yields
\begin{align*}
\lambda |T^\delta_j|/(4\log 1/\delta)^2&\le N^{-1}(2\log 1/\delta)^2\int_{T^\delta_j} \sum_{k=1}^{M_0}\chi_{\ii_k}(x) \, dx
\\
& =N^{-1}(2\log 1/\delta)^2\sum_{k=1}^{M_0} |T^\delta_{\ii_k}\cap T^\delta_j|
\le \frac{C\delta^3(\log 1/\delta)^2}{\theta N} \cdot M_0.
\end{align*}
Since $|T^\delta_j|\approx \delta^2$, this yields \eqref{4}, which finishes the proof of Proposition \ref{mainprop}.

\medskip

To finish matters, we still have to prove Lemmas \ref{auxest1} and \ref{bushlemma}.  Both require the following simple lemma which is essentially in \cite{S1}.  (See also \cite{MS}.)

\ble\label{inj}  Suppose that $\gamma_j$, $j=1,2$ are geodesics of length $\alpha\le \min \{1, (\text{inj }M^3)/2\}$ and assume that the $\gamma_j$ belong to a fixed compact subset $K$ of $M^3$.  Suppose also that $a\in T^\delta_{\gamma_1}\cap T^\delta_{\gamma_2}$.  Then there is a constant $c>0$, depending on $(M^3,g)$ and $K$, but not on $0<\lambda, \delta \le 1$, so that
$$(T^\delta_{\gamma_1}\cap T^\delta_{\gamma_2})\backslash B(a,\lambda)=\emptyset
\quad \text{if } \, \text{angle }(T^\delta_{\gamma_1},T^\delta_{\gamma_2})\ge \delta/c\lambda.$$
\ele

\noindent{\it Proof of Lemma \ref{bushlemma}.}  If we sum \eqref{bass2} we conclude that $\sum_{j=1}^{M_0}|E\cap T^\mu_j|\ge c_0M_0\mu^2\rho$, for some fixed $c_0>0$.  From this, we conclude that there must be a point $a\in E$ belonging to at least 
$$N_0=c_0M_0\mu^2\rho/|E|$$
of the tubes $T^\mu_j$.  Label these as $\{T^\mu_{j_k}\}_{1\le k\le N_0}$.

If we invoke the preceding lemma, we conclude that if $C_1$ in \eqref{bass1} is large enough then $(T^\mu_{j_k}\cap T^\mu_{j_{k'}})\backslash B(a,\mu/\theta)=\emptyset$ if $k\ne k'$.  Hence, by \eqref{bass2}
$$|E|\ge |E\cap(\cup (T^\mu_{j_k}\backslash B(a,\mu/\theta)))|\ge N_0\rho\mu^2=c_0M_0\mu^4\rho^2/|E|,$$
which of course yields \eqref{bconcl}.

\medskip

\noindent{\it Proof of Lemma \ref{auxest1}.}  We first notice that \eqref{good} would clearly follow from the dyadic estimates
\begin{multline}\label{c}
\bigl(\int_{\theta/2\le |x'|\le \theta}|\Adel f(x')|^2\, dx'\bigr)^{1/2}\le C(\log 1/\delta)^{1/2}\|f\|_{L^2}, 
\\ \delta\le \theta, \mu\le 1/2, \, \, \supp f\subset \{(y_1,y'): |y'|\in[\mu/2,\mu]\}.
\end{multline}
For this, write $\{x': |x'|\in [\theta/2,\theta]\}=\cup_{1\le k\le 10\theta/\delta}\Pi^\delta_k$, where 
$$\Pi^\delta_k=\{y\in {\Bbb R}^2: |<\nu_k,y>|\le \delta\}$$
for a collection of $\delta/\theta$-separated points $\nu_k\in S^1$.  We then have that
\begin{equation}\label{d}
\sum_{1\le k\le 10\theta/\delta}\chi_{\Pi^\delta_k}(x')\le C, \quad |x'|\in
[\theta/2,\theta],
\end{equation}
for some uniform constant $C$.

On the other hand, if $x'\in \Pi^\delta_k$ and $\gamma_{x'}\cap \gamma_0\ne \emptyset$, where $\gamma_0$ is the common geodesic,  then 
\begin{equation}\label{crucial}
T^\delta_{\gamma_{x'}}\subset \{y: |<y',\nu_k>|\le C\delta\},
\end{equation}
for some fixed constant $C$, if we are in the constant curvature case.  This just follows from the fact that every Fermi two plane must be totally geodesic under this assumption.  The next thing we need to use is that
\begin{equation}\label{e}
\sum_{1\le k\le 10\theta/\delta}\chi_{\Pi^{C\delta}_k}(y')\le C'\theta/\mu, \quad
\text{if } |y'|\in [\mu/2,\mu].
\end{equation}
Next, since the damping factors $w_{\gamma_{x'}}(y)$ in the definition of $\Adel$ are $\le C(\mu/\theta)^{1/2}$ when $|y'|\in [\mu/2,\mu]$, we conclude using \eqref{d} and \eqref{e} and a twofold application of Schwarz's inequality (and possibly replacing $\delta$ by a fixed multiple of $\delta$) that \eqref{c} would follow from showing that, when $\nu\in S^1$,
\begin{multline}\label{f}
\bigl(\int_{|<x',\nu>|\le \delta}|\Adel f(x')|^2\, dx'\bigr)^{1/2}\le C(\log 1/\delta)^{1/2}\|f\|_{L^2},
\\
\text{if } \, f(y)=0 \quad \text{when } \quad |<y',\nu>|\ge \delta.
\end{multline}

To prove this we need to appeal to a variable coefficient version of a theorem of C\'ordoba \cite{C} which is essentially in \cite{MSS}.   To state it we now suppose that $(M^2,g)$ is a two-dimensional Riemannian manifold.  If we fix a geodesic $\gamma_0\subset M^2$ of length $\alpha\le \min \{1, (\text{inj }M^2)/2\}$, we consider all geodesics $\{\gamma\}$ of this length which are close to $\gamma_0$.  If $\gamma_1(t)$ is another geodesic which intersects $\gamma_0$ orthogonally and is parameterized by arclength, we set
$$g_\delta^*(t)=\sup_{\gamma\ni\gamma_1(t)}\delta^{-1}\int_{\{y: \, \text{dist }(y,\gamma)\le \delta \}} |g(y)|\, dy.$$
Then the estimate we require is
\begin{equation}\label{g}
\|g_\delta^*\|_{L^2(dt)}\le C(\log1/\delta)^{1/2}\|g\|_{L^2(M^2)}, 
\end{equation}
assuming as usual that the functions involved are supported in a fixed compact set $K$.

The preceding estimate implies \eqref{f} if we use once more the fact that Fermi two planes are totally geodesic when $M^3$ has constant curvature.

To prove \eqref{g}, it suffices to prove a linearized version.  Specifically, if we choose for each $t$ a geodesic $\gamma_t$ as above which contains the point $\{\gamma_1(t)\}$, it suffices to show that the operator
$$(Tg)(t)=\delta^{-1}\int_{\{y\in K: \, \text{dist }(y,\gamma_t)\le \delta\}}g(y) \, dy$$
is bounded from $L^2(K)$ to $L^2(dt)$ with norm $\le C(\log 1/\delta)^{1/2}$, with $C$ being a uniform constant.  This in turn would follow if and only if we had bounds for $TT^*$ of the form
\begin{equation}\label{h}
\|TT^*g\|_{L^2(dt)}\le C(\log 1/\delta) \|g\|_{L^2(dt)}.
\end{equation}

To verify this assertion we note that $TT^*$ has kernel
$$K(t,s)=\delta^{-2}|T^\delta_s\cap T^\delta_t|$$
if $T^\delta_s=\{y\in K: \text{dist }(y,\gamma_t)\le \delta\}$.  Consequently, we have the trivial estimate that $K(t,s)=O(\delta^{-1})$.  If we use Lemma \ref{inj} we also get $K(t,s)=O(|t-s|^{-1})$ since if $T^\delta_s\cap T^\delta_t\ne \emptyset$, we must have that $\text{angle }(T^\delta_s,T^\delta_t)\ge c|t-s|$ for some uniform $c>0$.  Since Young's inequality and these two estimates for $K(s,t)$ give \eqref{h} the proof is complete.

\medskip

\noindent{\bf Remark.}  The only ingredient in the proof of Theorem 2.1 which used the constant curvature assumption was Lemma \ref{auxest1}.  The only step in the proof of the latter result which used our hypothesis was \eqref{crucial}.  What we really used was that if $\gamma_0\subset M^3$ is a geodesic of length $\alpha$ as above, then we can choose local coordinates near $\gamma_0$ so that every resulting two plane which contains $\gamma_0$ is (locally) totally geodesic.
Unfortunately, this is true for all such $\gamma_0$ if and only if $M^3$ has constant curvature (see Proposition \ref{taylor} below).  We shall see later that for variably curved manifolds Lemma \ref{auxest1} always breaks down.  Thus, bounds like \eqref{good} can only hold in the special cases where $M^3$ is isometric near every point to Euclidean space, a sphere or hyperbolic space.  Nonetheless, we shall be able to prove some new results concerning lower bounds on the dimension of Nikodym-type sets in a generic class of manifolds by proving a much weaker auxiliary maximal estimate.

\newsection{Chaotic curvature and spaces of variable curvature}

Recall that in three dimensions, a connected Riemannian manifold $M^3$ has constant curvature if and only if its Einstein tensor $B_{ij}$ vanishes identically.  This, we recall, is just the trace free part of the Ricci tensor, $R_{ij}$, that is, 
$$B_{ij}=R_{ij}-Rg_{ij}/3,$$
where $R=\sum R^i{}_{i}$ is the scalar curvature.\footnote{ Here we are using the usual convention of lowering and raising indices using the metric.} 
In many ways the Einstein tensor measures the amount of symmetry of $M^3$.

In the preceding section we saw that in the case of maximal symmetry where $B_{ij}\equiv 0$, we must have that $\text{dim }\Omega\ge 5/2$ if $\Omega$ is a Nikodym-type set.  On the other hand, as we noted in the introduction, such results are not stable under arbitrarily small perturbation.  Indeed, one can construct arbitrarily small perturbations of any constant curvature manifold so that, in the resulting manifold $M^3$, there are Nikodym-type sets with $\text{dim }\Omega=2$.

This suggests that if we wish to have any improvements over the lower bound $\text{dim }\Omega\ge 2$ in \cite{MS} for spaces of non-constant curvature, we need to impose a condition which will ensure that we stay away from the symmetric case.  Such a condition would have to involve the ``off-diagonal'' parts of $R_{ij}$ or $B_{ij}$.

Let us be more specific.  Suppose that $\gamma(t)$, $0\le t\le t_0$ is a geodesic parameterized by arclength.  Suppose further that for $0\le t\le t_0$ the vectors $X(t)\in T_{\gamma(t)}M^3$ are orthogonal to $\dot \gamma(t)$ and also parallel along $\gamma$.
(Recall that the latter means that
$$DX^i/\partial t = \partial X^i/\partial t+\sum \Gamma_{jk}{}^i X^kd\gamma^j/\partial t=0, \, \, \, \forall i,$$ 
if $\Gamma_{jk}{}^i$ are the Christoffel symbols of the metric.)
If then for $Y(t)\in T_{\gamma(t)}M^3$ we let $Y_\perp(t)$ be the projection onto the orthogonal compliment of the space spanned by $X(t)$ and $\dot \gamma(t)$,
we can place a condition on the ``off-diagonal'' part of $R_{ij}$ by requiring that
\begin{equation}\label{3.1}
|Z(t)|+|DZ/\partial t|\ne 0, \, \, \,
\text{if } \, \, Z(t)=Y_\perp (t)\, \, \text{with } \, \,
Y^i(t)=\sum R^i{}_j(\gamma(t))X^j(t).
\end{equation}
Also, since $W_\perp (t)=0$ if $W^i(t)=\sum g^i{}_j(\gamma(t))X^j(t)$, we get the same condition if we use the Einstein tensor instead of the Ricci tensor in \eqref{3.1}.

\bdeff\label{def3.1}  We shall say that $M^3$ satisfies the {\it chaotic curvature} condition at $x_0\in M^3$ in the direction $\tau_0\in T_{x_0}M^3\backslash 0$ if \eqref{3.1} holds  whenever $\gamma$ and $X$ are as above with $x_0\in \gamma$ and $\dot \gamma=\pm \tau_0/\|\tau_0\|$ at $x_0$. 
We shall say that a given family ${\cal F}$ of geodesics of a given length $\alpha$ satisfies the chaotic curvature condition if ${\cal F}$ is closed and \eqref{3.1} holds for every $\gamma\in {\cal F}$.  We shall say that $M^3$ is variably curved if this condition holds for all geodesics.
\edeff

Before moving on, we should explain how \eqref{3.1} simplifies if one works in a Fermi normal coordinate system about $\gamma$.  First of all, in such coordinates, the vector fields $X(t)$ above must be of the form $a_2\partial/\partial x_2 +a_3 \partial/\partial x_3$ with $a_j$ constant.  If $X(t)=\partial/\partial x_2$, then \eqref{3.1} just means that $R_{23}$ can only vanish to first order on $\gamma$.  That is, if $R_{23}(t)$ denotes the $23$ component of the Ricci tensor at $\gamma(t)$, then
\begin{equation}\label{3.2}
|R_{23}|+|dR_{23}/dt|\ne 0.
\end{equation}
We can replace this condition in the way that it will be used later if we recall that the Christoffel symbols vanish at the center of Fermi normal coordinates and hence, on $\gamma$,  $2R_{23}=g_{23,12}+g_{12,23}-g_{13,22}-g_{22,13}$, with
$g_{jk,i_1\dots i_m}=\partial/\partial x_{i_1}\cdots \partial/\partial x_{i_m}g_{jk}$.  If we now use \eqref{fermi2}, we conclude that \eqref{3.2} holds if and only if
$$|g_{11,23}|+|g_{11,231}|\ne 0 \quad \text{on } \, \, \gamma.
$$

Similar reasoning gives that if $X(t)=\cos \psi \partial/\partial x_2+\sin\psi \partial/\partial x_3$, then \eqref{3.1} becomes
\begin{multline}\label{3.3}
|(\cos\psi \partial/\partial x_2\! +\! \sin\psi \partial/\partial x_3)
(\sin\psi \partial/\partial x_2\! -\! \cos\psi \partial/\partial x_3) g_{11}|
\\
+
|(\cos\psi \partial/\partial x_2\! +\! \sin\psi \partial/\partial x_3)
(\sin\psi \partial/\partial x_2\!-\! \cos\psi \partial/\partial x_3) g_{11,1}|\ne 0, \, \, \, \text{on } \, \gamma.
\end{multline}
Clearly, if we fix $t_0$, then we can always choose $\psi$ so that the first term vanishes at $\gamma(t_0)$, which explains our condition involving the next best thing that these terms can only vanish to first order along $\gamma$.

There is another way of seeing that we can always choose $X(t)$ so that the first term in \eqref{3.1} vanishes at a given point $\gamma(t_0)$.  We first recall that if we work in a given Fermi normal coordinate system and choose $X(t)$ as in the preceding step, then at $\gamma(t_0)$ twice the first term in \eqref{3.1} is the sectional curvature for the two-plane spanned by $\partial/\partial x_1$ and $\cos(\psi+\pi/4)\partial/\partial x_2+\sin(\psi+\pi/4)\partial/\partial x_3$ minus the sectional curvature for the two-plane spanned by $\partial/\partial x_1$ and $\cos(\psi-\pi/4)\partial/\partial x_2+\sin(\psi-\pi/4)\partial/\partial x_3$.  Since this difference is a function of $\psi$ which clearly has mean value zero, we can always choose $\psi$ so that it vanishes.  This formulation explains our choice of the phrase ``variably curved'' in Definition 3.1.

Before stating the main result of this section, let us explain how the condition \eqref{3.1} arises naturally in the applications we have in mind.  To do this we recall that the counterexamples in \cite{MS} showing that one can have unfavorable lower bounds for Nikodym-type sets all involved a family of space-filling geodesics which become highly focused in a lower dimensional submanifold ${\cal H}$.  The following result shows that for variably curved manifolds there can only be third or fourth order focusing if ${\cal H}$ is a Fermi two-plane, and, hence, in contrast to the constant curvature case, Fermi two-planes are of course not totally geodesic in this case.

\bprop\label{taylor}  Fix a local Fermi normal coordinate system about a given geodesic $\gamma_0\subset M^3$.  We then fix $-\pi\le \theta<\pi$ and small $x_1$ and let $\gamma=\gamma_{x_1,\theta}$ be the geodesic parameterized by arclength satisfying
$$\gamma(0)=(x_1,0,0), \quad \text{and } \, \, \dot \gamma(0)=(\cos\theta,\sin\theta,0).$$
Then,
\begin{equation}\label{3.4}
d^2\gamma(0)/dt^2=0,
\end{equation}
and, moreover, if $\gamma^3$ denotes the third coordinate of $\gamma$,
\begin{align}\label{3.5}
d^3\gamma^3(0)/dt^3&=-\frac12 \cos^2\theta\sin\theta g_{11,23}(x_1,0,0) -3\cos\theta\sin^2\theta g_{12,23}(x_1,0,0) 
\\
\label{3.6}
d^4\gamma^3(0)/dt^4&=-\cos^3\theta\sin\theta g_{11,123}(x_1,0,0)+O(\theta^2).
\end{align} 
\eprop

\begin{pf}  We shall use the classical Jacobi equation
\begin{equation}\label{3.7}
d^2\gamma^k/dt^2=\sum_{i,j}\Gamma_{ij}{}^k \dot\gamma^{i}\dot\gamma^j,
\end{equation}
where the Christoffel symbols are evaluated at $\gamma$ and are defined by
$$\Gamma_{ij}{}^k=\sum_l \Gamma_{ijl}g^{lk}, \quad
2\Gamma_{ijk}=g_{ik,j}+g_{jk,i}-g_{ij,k}.$$
If we use \eqref{fermi2} we conclude that
\begin{equation}\label{3.8}
(\partial/\partial x_2)^m(\partial/\partial x_3)^n(\Gamma_{ijk}-\Gamma_{ij}{}^k)=0 \, \, \text{if }x_2=x_3=0, \, \, \text{and } \, m+n\le 2,
\end{equation}
and also that $\Gamma_{ij}{}^k(x_1,0,0)=0$.  The latter and \eqref{3.7} yield \eqref{3.4}.   If we use this and \eqref{3.8} we conclude that
\begin{equation}\label{3.9}
d^3\gamma^3(0)/dt^3=\sum_{j,k,l}\Gamma_{jk3,l}(x_1,0,0)
\dot\gamma^j(0)\dot\gamma^k(0)\dot\gamma^l(0),
\end{equation}
and
\begin{equation}\label{3.10}
d^4\gamma^3(0)/dt^4=\sum_{j,k,l,m}\Gamma_{jk3,lm}(x_1,0,0)
\dot\gamma^j(0)\dot\gamma^k(0)\dot\gamma^l(0)\dot\gamma^m(0).
\end{equation}
Because of the initial conditions, in both cases the summands where one of the indices is $3$ vanish.  
Also, clearly $\Gamma_{jk3,1}(x_1,0,0)=0$, and using \eqref{fermi1} one finds that $\Gamma_{223,2}(x_1,0,0)=0$.  Therefore, since $\Gamma_{123}=\Gamma_{213}$
$$d^3\gamma^3(0)/dt^3=\cos^2\theta\sin\theta \Gamma_{113,2}(x_1,0,0)+2\cos\theta \sin^2\theta \Gamma_{123,2}(x_1,0,0).$$
If we use \eqref{fermi2} we conclude that $g_{jk,m1}=0$ and hence $\Gamma_{113,2}=-g_{11,23}/2$ at $(x_1,0,0)$.  For the other term, we need to use \eqref{fermi1} to conclude that $g_{13,22}=-2g_{12,23}$ and consequently $\Gamma_{123,2}=-3g_{12,23}$ at $(x_1,0,0)$.  By combining these calculations we get \eqref{3.5}.

The proof of \eqref{3.6} is similar.  Modulo $O(\theta^2)$ terms, the sum in the right side of \eqref{3.10} just involves terms where one of the $jklm$ is $2$ and the rest are $1$.  Thus,
$$d^4\gamma^3(0)/dt^4=(2\Gamma_{113,12}+2\Gamma_{123,11})\cos^3\theta\sin\theta +O(\theta^2),$$
where the Christoffel terms are evaluated at $(x_1,0,0)$.  Since $\Gamma_{123,11}=0$ and $2\Gamma_{113,12}=-g_{11,123}$ at this point, we obtain \eqref{3.6}.
\end{pf} 

Before moving on we should note for later use that we can use this result and a simple rotation argument to compute the Taylor coefficients of geodesics which are initially tangent to other Fermi two-planes.  Specifically, we have the following

\bcor\label{cor3.3}
Assume that coordinates are chosen as above and let $\gamma_{x_1\theta\psi}$ be the geodesic parameterized by arclength satisfying
$$\gamma_{x_1\theta\psi}(0)=(x_1,0,0), \quad
\text{and } \, \, \dot\gamma_{x_1\theta\psi}(0)=(\cos\theta,\cos\psi\sin\theta,\sin\psi\sin\theta).$$
Then
\begin{equation}\label{3.11}
d^2\gamma_{x_1\theta\psi}(0)/dt^2=0,
\end{equation}
and moreover, if $\gamma^\perp_{x_1\theta\psi}=<\gamma_{x_1\theta\psi}, (0,-\sin\psi,\cos\psi)>$, then
\begin{align}\label{3.12}
d^3\gamma^\perp_{x_1\theta\psi}(0)/dt^3&=\rho(x_1,\psi)
\cos^2\theta\sin\theta +O(\theta^2),
\\
\label{3.13}
d^4\gamma^\perp_{x_1\theta\psi}(0)/dt^4&=2\rho'_{x_1}(x_1,\psi)
\cos^3\theta\sin\theta 
+O(\theta^2),
\end{align}
where 
\begin{equation}
\label{3.14}
2\rho(x_1,\psi)
=(\cos\psi \partial/\partial x_2\! +\! \sin\psi \partial/\partial x_3)
(\sin\psi \partial/\partial x_2\! - \! \cos\psi \partial/\partial x_3) g_{11}(x_1,0,0).
\end{equation}
\ecor

Note that if $M^3$ is variably curved then for any fixed $\psi$ the function $x_1\to \rho(x_1,\psi)$ can only have first order zeros.

Having gone through the preliminaries we now state our main result.

\bth\label{thm3.4}
Assume that $M^3$ is variably curved (see Definition 3.1).  Then if $\Omega\subset M^3$ is a Nikodym-type set in the sense of Definition 1.1 its Minkowski dimension satisfies
\begin{equation}\label{3.15}
\text{dim }\Omega \ge 7/3.
\end{equation}
More generally, if ${\cal F}$ is a family of geodesics of length $\alpha$ satisfying the chaotic curvature condition, then the same conclusion holds if for $\lambda$ sufficiently close to $1$
$$|\{x\in M^3: \, \exists \gamma_x\in {\cal F} \, \, \text{with } \, \,
x\in \gamma_x \, \, \text{and } \, \, |\gamma_x\cap \Omega|\ge \lambda |\gamma_x|\}|>0.$$
\eth

To prove this we need suitable estimates for the associated Nikodym maximal 
operator:
\begin{equation}\label{nikmax2}
f^*_\delta(x)=\sup_{x\in\gamma_x\in{\cal F}} 
|T_{\gamma_x}^\delta|^{-1} \int_{T_{\gamma_x}^\delta}|f(y)|\, dy.
\end{equation}
If
we argue as in the proof of Corollary 2.2, we then see that \eqref{3.15} is a consequence of the following 

\bprop\label{prop3.5}  Let $E$ be contained in a compact subset of a coordinate patch.  Then if $\varepsilon>0$ there is a uniform constant $C_\varepsilon$ so that, for all $\lambda$ sufficiently close to $1$,
\begin{equation}\label{3.16}
\delta^{2/3+\varepsilon}|\{x: \, f^*_\delta(x)\ge \lambda\}|\le C_\varepsilon |E|, \quad f=\chi_E.
\end{equation}
\eprop

The fact that, at present we can only prove estimates like \eqref{3.16} with $\lambda$ close to $1$ accounts for the same requirement in our definition of Nikodym-type sets.  It is also reflected in the fact that the lower bound \eqref{3.15} involves the Minkowski dimension, rather than the Hausdorff dimension.  In the next section, though, we shall see that we can drop this assumption if in the left we replace $\delta^{2/3}$ by $\delta^{3/4}$ and allow constants with an unfavorable dependence on $\lambda$.

The reason for this limitation is that we can only prove rather weak estimates for the auxiliary maximal function arising in the proof.  Specifically, let us fix a ``common geodesic'' $\gamma_0$ of length $\alpha$, and, as before, choose Fermi normal coordinates about it so that 
$$\gamma_0=\{(s,0,0): 0\le s\le \alpha\}.$$
We then modify \eqref{auxdef} by setting\footnote{Note that since we are not proving estimates for small scales $\lambda$, the weights $w_{\gamma_{x'}}$ are no longer relevant.  As we shall see at the end they are needed if one wishes to prove what seem to be the optimal estimates.}
\begin{equation}\label{3.17}
(\Adel f)(x')=
\sup_{\{\gamma_{x'}\ni (0,x'): \, \gamma_{x'}\cap\{(s,0,0): \alpha/2\le s\le \alpha\}\ne \emptyset \}}
|T^\delta_{\gamma_{x'}}|^{-1}\int_{T^\delta_{\gamma_{x'}}} 
|f(y)|\, dy.
\end{equation}
The main estimate on which \eqref{3.15} and \eqref{3.16} then is contained in the following

\bprop\label{prop3.6}  Suppose that $\gamma_0$ satisfies the chaotic curvature condition.  Then there is an $r>0$ so that if $\varepsilon>0$ and if $\lambda$ is sufficiently close to $1$ then
\begin{equation}\label{3.18}
|\{x'\in B(0,r):\, (\Adel f)(x')\ge \lambda \}|\le C_\varepsilon \delta^{-1/3-\varepsilon}|E|, \quad f=\chi_E,
\end{equation}
for some uniform constant $C_\varepsilon$.
\eprop

The restriction $|x'|\le r$ in \eqref{3.18} forces the averages in \eqref{3.17} to just involve tubes whose centers  $\gamma_{x'}$ are close to $\gamma_0$.

Clearly the proof of \eqref{weak} can be adapted to show that \eqref{3.18} implies \eqref{3.16}.  Let us sketch the argument.  In the present context one must show that
\begin{equation}\label{3.19}
M\delta^2\le C_\varepsilon (\delta^{-2/3-\varepsilon}|E|)^{4/3},
\end{equation}
assuming that \eqref{assumption} holds with $\lambda$ close to $1$.  As before we may assume that the tubes $T^\delta_j$ are close to one another, which allows us to use Proposition \ref{3.6}.

For the next step one needs to modify Lemma 2.4 slightly.  Using the pigeonhole principle just as before we argue that if \eqref{assumption} holds then we can find $N\in {\Bbb N}$ and small $\theta>0$ so that \eqref{case1} holds for $M/2$ values of $j$, while now
$$|\{x\in T^\delta_j\cap E: \, \text{card }{\cal I}_\theta(x,j)\ge N/2\}|\ge (4\log_21/\delta)^{-1}\lambda |T^\delta_j|,$$
for at least $M/(2\log_2 1/\delta)$ indices $j$.  As we noted before this sort of result does not rely on the geometry of $M^3$.  For similar reasons \eqref{trivial} remains valid.

The remaining key estimate in our proof, however, does rely on our geometric assumptions.  Specifically, the analog of \eqref{nontrivial} which we can prove, says that, for small $\delta$,
\begin{equation}\label{3.20}
|(E\backslash B(a,\delta^\varepsilon))\cap T^\mu_j|\ge C_\varepsilon \mu \delta^{1+1/3+\varepsilon}N, \quad a\in M^3, \quad \varepsilon>0.
\end{equation}
One can then use the bush lemma, Lemma 2.6, just as before to deduce \eqref{3.19}.

To prove \eqref{3.20} one first argues as before that it suffices to prove the lower bound when $E\backslash B(a,\delta^\varepsilon)$ is replaced by $E$.  One can then easily adapt the proof of \eqref{nontrivial} to verify the resulting inequality.  One just notices that if we work in Fermi normal coordinates about $\gamma_j$ then our assumptions imply lower bounds for 
${\cal A}_{C\delta} f(x'_i)$, assuming that the tube $T^\delta_i$ intersects the top half of $T^\delta_j$.  Specifically, if we assume that \eqref{a} holds and replace \eqref{0} by $|T^\delta_i \cap E|\ge \rho|E|$ with $\rho=\lambda$, then the lower bound in 
\eqref{pig} is just replaced by ${\cal A}_{C\delta} f(x'_i)\ge c\rho$. (One does not divide by $\log 1/\delta$ now since the multiplicity argument involves ${\cal I}_\theta(x,j)$, rather than ${\cal I}_{\theta\mu}(x,j)$ as in \eqref{case2}.)
Since the discrete version of \eqref{3.18} then yields the aforementioned variant of \eqref{3.20}, the proof is complete.

\medskip

Let us conclude this section by proving Proposition \ref{prop3.6}.  We shall obtain \eqref{3.17} from an $L^2$-estimate involving a somewhat weaker maximal operator.  Specifically, let us set
\begin{equation}\label{3.21}
({\cal W}_\delta f)(x')=
\sup_{\{\gamma_{x'}\ni (0,x'): \, \gamma_{x'}\cap\{(s,0,0): \alpha/2\le s\le \alpha\}\ne \emptyset \}}
|T^\delta_{\gamma_{x'}}|^{-1}\int_{T^\delta_{\gamma_{x'}}} a_{\gamma_{x'}}(y)
|f(y)|\, dy,
\end{equation}
where we assume that the weights satisfy $0\le a_{\gamma_{x'}}\le c_0^{-1}$ 
and
\begin{equation}\label{3.22}
|\{y\in T^\delta_{\gamma_{x'}}: \, a_{\gamma_{x'}}(y)\ge 1\}|\ge c_0|T^\delta_{\gamma_{x'}}|, \quad \varepsilon>0,
\end{equation}
for some uniform constant
$$0< c_0\le 1.$$
We then claim that estimates of the form
\begin{equation}\label{3.23}
\|{\cal W}_\delta f\|_{L^2(|x'|\le r)}\le C_\varepsilon \delta^{-1/6-\varepsilon}\|f\|_{L^2}
\end{equation}
imply \eqref{3.18} for $\lambda$ close to $1$.

This is easy to check.  We first note that \eqref{3.23} of course implies that the weak-type bounds
$$|\{x'\in B(0,r): \, {\cal W}_\delta f(x')\ge \mu\}|\le C_\varepsilon \delta^{-1/3-2\varepsilon}\mu^{-2}\|f\|_{L^2}^2, \quad \varepsilon >0.$$
But then this yields \eqref{3.18} for $\lambda$ close to $1$ since
$$\{x'\in B(0,r): \, \Adel f\chi_{E}(x')\ge \lambda\}\subset
\{x'\in B(0,r): \, {\cal W}_\delta \chi_E(x')\ge c_0-(1-\lambda)\}.$$

Estimate \eqref{3.23} is somewhat similar to an estimate of Mockenhaupt, Seeger and the author \cite{MSS} that involved a maximal function arising from averages over tubes around null geodesics in a curved Lorentz manifold.  Even though \eqref{3.23} seems more complicated than the auxiliary estimate in \cite{MSS} we can follow the general strategy in that paper of proving our bounds using a simple square function argument along with $L^2$ estimates for the underlying Fourier integral operators.  In \cite{MSS}, after a change of variables, we could use an $L^2$ estimate of H\"ormander \cite{H1} for Fourier integral operators whose canonical relation is locally a canonical graph.  Here, though, we need to appeal to an $L^2$ estimate of Melrose and Taylor \cite{MT} which says that there are $L^2$ bounds with a loss of $1/6$ derivatives for Fourier integral operators with two-sided folds.  The loss of $1/6$ derivatives accounts for the loss of $1/6$ on the lower bounds for the dimension of Nikodym-type sets in Theorem \ref{thm3.4} versus the corresponding bound (2.5) for the constant curvature case.  To be able to apply the results of Melrose and Taylor we shall use Corollary \ref{cor3.3}, which concerns the geometry of the curves $\gamma_{x'}$.  The fact that the function $\rho$ there may vanish complicates matters and this is the main reason that we, for instance, can only prove bounds for certain operators ${\cal W}_\delta$ rather than for the more natural operators of the form \eqref{auxdef}.

The Fourier integral estimates will allow us to handle ${\cal W}_\delta f(x')$ when $|x'|$ is not to small.  Specifically, we want to avoid the trivial region where $|x'|\approx \delta$.  However, since the bounds in \eqref{3.23} involve $\delta^{-1/6}$ we can handle very small $|x'|$ using trivial arguments.  We just notice that $T^\delta_{\gamma_{x'}}\subset \{y=(y_1,y'): \, |y'|\le C|x'| \}$ for some uniform constant $C$.  Thus, using Schwarz's inequality, one can check that $({\cal W}_\delta f)(x')\le C\delta^{-2}|x'| \, \|f\|_2.$  As a result,
$$\|{\cal W}_\delta f\|_{L^2(|x'|\le \delta^{11/12})}\le C\delta^{-1/6}\|f\|_{L^2}.$$

We need another easy reduction.  This is needed since we shall want to scale the operators based on the size of $|x'|$.  To allow this, we notice that it suffices to show that we have uniform estimates over dyadic annuli.  Specifically, it is enough to show that
\begin{equation}\label{3.24}
\|{\cal W}_\delta f\|_{L^2(2^{-k}\le |x'|\le 2^{-k+1})}\le C_\varepsilon \delta^{-1/6-\varepsilon }\|f\|_2, \quad \varepsilon >0, \, \, \delta^{11/12}\le 2^{-k}\le r.
\end{equation}
Having $\delta^{11/12}$ here is not so important; for the arguments to follow it is just convenient to only have to prove the estimate for $|x'|\ge \delta^\sigma$ for some $\sigma<1$.

We now turn to the scaling argument.  Working in our Fermi normal coordinate system we shall want to scale $x'$ and the last two variables of $y$, $y'=(y_2,y_3)$ while keeping $y_1$ fixed.  To this end, let
\begin{equation}\label{3.25}
\gamma^k_{x'}=\{y: \, (y_1,2^{-k}y')\in \gamma_{2^{-k}x'}\},
\end{equation}
and
\begin{equation}\label{3.26}
a_{\gamma^k_{x'}}(y)=a_{\gamma_{x'}}(y_1,2^{-k}y').
\end{equation}
If we let $T^\mu_{\gamma^k_{x'}}$ be a $\mu$-neighborhood around $\gamma^k_{x'}$, put
\begin{equation}\label{3.27}
({\cal W}_{k,\mu}f)(x')=\sup_{\{\gamma^k_{x'}\ni (0,x'): \, \gamma^k_{x'}\cap \gamma_0\ne \emptyset \}} |T^\mu_{\gamma^k_{x'}}|^{-1}\int_{T_{\gamma^k_{x'}}} a_{\gamma^k_{x'}}(y) |f(y)|\, dy.
\end{equation}
Taking $\mu\approx 2^k\delta$, we conclude that \eqref{3.24} would follow from uniform bounds of the form
\begin{equation}\label{3.28}
\|{\cal W}_{k,\mu}f\|_{L^2(|x'|\in [1,2])}\le C_\varepsilon \mu^{-1/6-\varepsilon }\|f\|_2, \, \, \varepsilon >0, \, \, 0<\mu<1, \, \, \delta^{11/12}\le 2^{-k}\le r.
\end{equation}

Notice that for fixed $x'$  and $k$, the supremum in \eqref{3.27} involves the one-parameter family of curves $\{\gamma^k_{x'}\}$ satisfying $(0,x')\in \gamma^k_{x'}$ and $\gamma^k_{x'}\cap \gamma_0\ne \emptyset$.  Until now it has been convenient to suppress this extra parameter from the notation.  Before taking in into account, though, we should note that Corollary \ref{cor3.3} and its proof (see also \eqref{3.3.38} below) show that the one-parameter families $\{\gamma^k_{x'}\}$ actually tend to a limit as $k\to +\infty$.  Indeed, up to $O(2^{-k})$ error terms they agree with a family which is independent of $k$.  With this in mind, write the one-parameter family as $\{\gamma^k_{x',s}\}$, $0\le s\le 1$, where
\begin{equation}\label{3.29}
\gamma^k_{x',s}=\{(y_1,\Xi^k(x',s,y_1))\}
\end{equation}
where $\Xi^k$ takes its values in ${\Bbb R}^2$ and satisfies $|D^m_{x',s,y_1}\Xi^k|\le C_m$ if $|x'|\approx 1$.

To set up the square function argument fix $b\in C^\infty_0({\Bbb R}^2)$ satisfying $\Hat b\ge 0$ and $\hat b(s)\ge1$, $|s|\le 1$, where $\Hat b$ denotes the Fourier transform.  We then set
$$(W_{k,\mu}f)(x',s)=(2\pi)^{-1}\int_{{\Bbb R}^3}\int_{{\Bbb R}^2} e^{i<(y'-\Xi^k(x',s,y_1)),\xi>} a_{\gamma^k_{x',s}}(y)b(\mu\xi)f(y)\, d\xi dy.
$$
It then follows that $\sup_{0\le s\le 1}W_{k,\mu}f(x',s)$ dominates ${\cal W}_{k,\mu}f(x')$ if $f$ is nonnegative.  Consequently, it suffices to show that we can choose weights $a_{\gamma_{x'}}$ so that the resulting operators $W_{k,\mu}$ satisfy the uniform bounds
$$\|\sup_{0\le s\le 1}|W_{k,\mu}f(\, \cdot\, ,s)|\, \|_{L^2(|x'|\in [1,2])} \le C_\varepsilon \mu^{-\varepsilon }\|f\|_2, \, \, \varepsilon >0, \, \, 0<\mu\le 1, \, \, 2^{-k}\le r.$$

To proceed, we need to make one last dyadic decomposition.  For this, let us fix $\beta\in C^\infty_0([1/4,2])$ satisfying $\sum_{-\infty}^\infty \beta(2^jt)=1$, $t>0$.  If we then set
$$b^j_{k,\mu}(s,y,\xi)=a_{\gamma^k_{x',s}}(y)b(\mu\xi)\beta(2^{-j}|\xi|)$$
and
\begin{equation}\label{3.30}
(W^j_{k,\mu}f)(x',s)=(2\pi)^{-2}\iint e^{i<(y'-\Xi^k(x',s,y_1)),\xi>}b^j_{k,\mu}(s,y,\xi)f(y)\, d\xi dy,
\end{equation}
we claim that it suffices to show that
\begin{equation}\label{3.31}
\|\sup_{0\le s\le 1}|W^j_{k,\mu}f(\, \cdot \, ,s)|\, \|_{L^2(|x'|\in [1,2])}\le C2^{j/6}\|f\|_2, \quad 2^{-k}\le r.
\end{equation}

This implies the proceeding inequality since $W^j_{k,\mu}=0$ if $2^j$ is larger than a fixed multiple of $\mu^{-1}$ and since $W_{k,\mu}-\sum_{j\ge 0}W^j_{k,\mu}$ has a bounded kernel and hence the resulting maximal operator is bounded between any $L^p$ spaces.

If $(s,y)\to a_{\gamma^k_{x',s}}$ belongs to a bounded subset of $C^\infty$, then $b^j_{k,\mu}$ belong to a bounded subset of zero-order symbols supported in the region where $|\xi|\in [2^{j-2},2^{j+2}]$.  To exploit this we need to use the simple fact (see, e.g., \cite{Soggebook}, p. 75) that if $F\in C^1({\Bbb R})$ then
$$\sup_\lambda |F(\lambda)|^2\le |F(0)|^2+2(\int |F|^2 d\lambda)^{1/2}(\int |F'|^2d\lambda)^{1/2}.$$
Because of this, we would get \eqref{3.31} if we could show that
\begin{multline}\label{3.32}
\bigl(\int_0^1\int_{|x'|\in[1,2]}|(\partial/\partial s)^mW^j_{k,\mu}f(x',s)|^2\, dx' ds\bigr)^{1/2}\le C2^{j/6+(m-1)j/2}\|f\|_2, \\ m=0,1, \, \, 2^{-k}\le r.
\end{multline}

The operators in \eqref{3.30} are just dilates of an original operator, where, like in \eqref{3.25}, the prime variables, $x'$ and $y'$ are dilated by $2^k$ while the other variables $y_1$ and $s$ remain fixed.  We shall want to show that after applying this change of scale the resulting operators in \eqref{3.32} belong to a bounded class of Fourier integral operators of order $(m-1/2)$ with two-sided folding canonical relations.  If we could do this, then the remaining estimate, \eqref{3.32}, would follow from the theorem of Melrose and Taylor \cite{MT} concerning the $L^2$ mapping properties of such operators.  Indeed, if the above claims were verified one would simply use the fact that 
$2^{-j/6-(m-1)j/2} (\partial/\partial s)^m W^j_{k,\mu}$ belong to a bounded subset of Fourier integral operators of order $-1/6$ with two-sided folding canonical relations 
${\cal C}_k$.  Since the resulting arguments immediately give that ${\cal C}_k \subset T^*{\Bbb R}^3\backslash 0 \times T^*{\Bbb R}^3\backslash 0$ and that the lifted canonical one forms for the two factors do not vanish on ${\cal C}_k$, the preceding family of operators must be uniformly bounded on $L^2$ by \cite{MT}.

The main step of course will be to try to compute the canonical relations ${\cal C}_k$ of these operators and to verify that the associated left and right projections have folding singularities. Fortunately, these relations tend to a limiting relation as $k\to +\infty$.  To see this and to help us verify the other claims, it turns out to be instructive to compute the 
projections for the canonical relation associated with the undilated operators.  This would be parameterized by the phase function
$$
(y'-\Xi(x',s,y_1))\cdot \xi,$$
if $\Xi=\Xi^0$, so that the undilated geodesics \begin{equation}\label{3.33}\gamma_{x',s}=\gamma^0_{x',s}=\{(y_1,\Xi(x',s,y_1))\}
\end{equation}
are just those arising in the definition of ${\cal W}_\delta$.

Using the above phase function, we can write the associated canonical relation as
$${\cal C}=\{\, (s,x',\Xi'_s\cdot\xi,\Xi'_{x_2}\cdot\xi,\Xi'_{x_3}\cdot\xi, \,
y_1,\Xi,\Xi'_{y_1}\cdot\xi, \xi)\, \}.$$
Thus, our initial task will be to show that we can construct weights as in \eqref{3.22} so that on the supports of the symbols, the maps
\begin{align}\label{3.34}
\Pi_r(x',s,y_1,\xi)&=(y_1,\Xi,\Xi'_{y_1}\cdot\xi,\xi),
\\
\label{3.35}
\Pi_l(s,x', y_1,\xi)&= (s,x',\Xi'_s\cdot\xi,\Xi'_{x_2}\cdot\xi,\Xi'_{x_3}\cdot\xi)
\end{align}
have at most folding singularities.  Note that $\Pi_l$ and $\Pi_r$ are the projections of ${\cal C}$ onto the first and last six variables, respectively.  In proving our assertions regarding these maps we are allowed to change variables since maps with at most folding singularities are preserved under nondegenerate changes of coordinates.
The proof will also show that $\partial\Xi/\partial(s,x')$ has full rank, which implies the necessary technical facts that ${\cal C}\subset 
T^*{\Bbb R}^3\backslash 0 \times T^*{\Bbb R}^3\backslash 0$ and that the lifted canonical one forms do not vanish on ${\cal C}$.

Note that 
$$\Sigma =\{\gamma_{x',s}\}$$
is a 3-dimensional submanifold of the 4-dimensional manifold of all geodesics in $M^3$. Until now we have just been using the fact that each element of $\Sigma$ can be parameterized by its intersection with the hypersurface $\{(0,x')\}$ along with an extra parameter $s$, which can be taken to be the ``height'' of its intersection with the common geodesic $\gamma^0$.  However, if we wish to establish our claims regarding the canonical relation
$${\cal C}\subset T^*\Sigma\backslash 0\times T^*M^3\backslash 0,$$
it seems much more natural to make use of a different parameterization of $\Sigma$ which takes into account the value of $(\gamma,\dot\gamma)$ at the intersection of a given geodesic $\gamma\in\Sigma$ with the common geodesic $\gamma_0$.  By doing this we can hope to make use of our earlier elementary calculations in Corollary \ref{3.3}.

Let us be more specific.  We first note that Corollary \ref{cor3.3} implies that the geodesic satisfying $(x_1,0,0)\in \gamma$ and $\dot \gamma/|\dot \gamma|=\pm (\cos\theta,\cos\psi\sin\theta,\sin\psi\sin\theta)$ is of the form
\begin{multline}\label{3.36}
t\to \bigl(x_1+t\cos\theta,t\cos\psi\sin\theta -\sin\psi\tan\theta p(x_1,\psi;t\cos\theta),
\\
t\sin\psi \sin\theta +\cos\psi\tan\theta p(x_1,\psi; t\cos\theta)\bigr) 
+r(x_1,\psi,\theta;t),
\end{multline}
where if $\rho$ is as in \eqref{3.14}
$$p(x_1,\psi;\tau)=-\frac1{12}\rho(x_1,\psi)\tau^3-\frac1{24}\rho'_{x_1}(x_1,\psi)\tau^4$$
and where the first coordinate of the remainder term vanishes and also
$$r=O(\theta t^3)\quad \text{and}\quad
|<r, (0,-\sin\psi,\cos\psi)>|=O(\theta t^5)+O(\theta^2t^3).$$
To write this geodesic in the form \eqref{3.33} there is a natural change of variables.  One first of course sets
$$y_1=x_1+t\cos\theta.$$
If one then lets
$$(x_2,x_3)=\tan\theta(\cos\psi,\sin\psi),$$
the above geodesic \eqref{3.36} is of the form
$$(y_1,\Xi(x_1,x_2,x_3,y_1)),$$
where
\begin{equation}\label{3.3.38}
\Xi=\Xi_0+O(|(x_2,x_3)|^2(x_1-y_1)^3)
\end{equation}
with
\begin{multline}\label{3.37}
\Xi_0=
\\
\bigl( x_2(y_1-x_1+r_3)-x_3(p(x_1,\psi;y_1-x_1)+r_5), \, x_3(y_1-x_1+r_3) +x_2(p(x_1,\psi;,y_1-x_1)+r_5)\bigr)
\end{multline}
for remainders $r_j$, $j=3,5$ satisfying
$$\bigl|D^m_{x_2,x_3}D^n_{x_1,x_2}r_j\bigr|=O((x_1-y_1)^{j-|n|}),
\quad 0\le j\le n.$$
For this estimate to be valid we need to assume that $|(x_2,x_3)|$ and hence $\theta$ are bounded away from zero.

To simplify the calculations let us first compute the analogs, $\Pi^0_r$ and $\Pi^0_l$, of \eqref{3.34} and \eqref{3.35} where $\Xi$ is replaced by $\Xi_0$.  We shall show that these two maps have at most folding singularities above all but isolated points of the geodesic in \eqref{3.36}.  We shall then indicate how this calculation leads to the same result for $\Pi_r$ and $\Pi_l$ if $\theta$ is sufficiently small.  Note that we only need to consider small $\theta$ if, as above, the norm in \eqref{3.23} is taken over a small ball.

Turning to the calculations at hand, we note that in our coordinates
$$\Pi^0_r(x,y_1,\xi)=(y_1,\, \Xi_0, \,(1+\frac{\partial r_3}{\partial y_1})(x_2,x_3)\cdot\xi +(\frac{\partial p}{\partial y_1}+\frac{\partial r_5}{\partial y_1})(-x_3,x_2)\cdot \xi,\, \xi),$$
and
\begin{multline*}
\Pi^0_l=
\\
(x,-(x_2,x_3)\cdot\xi-\frac{\partial p}{\partial x_1}(-x_3,x_2)\cdot \xi, (y_1-x_1+r_3)\xi_1+(p+r_5)\xi_2,(y_1-x_1+r_3)\xi_2-(p+r_5)\xi_1)
\\
+ (0,0,0, \frac{\partial r_3}{\partial x_1},\frac{\partial r_3}{\partial x_2},\frac{\partial r_3}{\partial x_3})(x_2,x_3)\cdot\xi
+(0,0,0, \frac{\partial r_5}{\partial x_1},\frac{\partial r_5}{\partial x_2},\frac{\partial r_5}{\partial x_3})(-x_3,x_2)\cdot\xi
.
\end{multline*}
We are abusing the notation somewhat when we for instance write $\frac{\partial p}{\partial y_1}$ instead of $\frac\partial{\partial y_1}p(x_1,\psi;y_1-x_1)$.

To proceed, we should recall the definition of a fold.  Specifically, if $\chi: {\Bbb R}^d\to {\Bbb R}^d$ is a smooth map, then $\chi$ is said to have a folding singularity at $x_0$ if $\text{rank }\chi'=d-1$, where $\chi'$ denotes the Jacobian, and moreover
\begin{equation}\label{3.38}
\text{Hess }\chi(x_0)=\bigl|\sum_{1\le j,k\le d}X_j X_k\frac{\partial^2}{\partial x_j\partial x_k} \, <\chi,Y>\bigr|\ne 0
\quad \text{at }\, x_0,
\end{equation}
if
$$X=(X_1,\dots,X_d)\in \text{Ker }\chi'(x_0)\cap S^{d-1} \, \, \,
\text{and } \, \, Y\in \text{Ker }(\chi'(x_0))^t\cap S^{d-1}.$$
The condition \eqref{3.38} is on the Hessian of $\chi$ at $x_0$ which should be thought of as a map from the kernel of $\chi'$ to the cokernel of $\chi'$.  It is well known and not hard to check that if $\chi$ has at most folding singularities then the same is true for its pullback under any diffeomorphism.  Consequently we are allowed to change coordinates if we wish to show that $\Pi^0_l$ or $\Pi^0_r$ has at most folding singularities at a given point.

To handle $\Pi^0_r$ we note that, in view of the form of the first and last two variables, $\Pi^0_r$ has at most a fold at $(x,y_1,\xi)$, $\xi\ne0$ if and only if when this $y_1$ and $\xi$ are fixed the map
$$x\to (\Xi_0,(1+\frac{\partial r_3}{\partial y_1})(x_2,x_3)\cdot\xi+(\frac{\partial p}{\partial y_1}+\frac{\partial r_5}{\partial y_1})(-x_3,x_2)\cdot \xi)$$
has at most a folding singularity at the above $x$.
Recalling the form of $\Xi_0$ in \eqref{3.37} we see that it is convenient to change variables by letting
$$(z_1,z_2,z_3)=(x_1,(y_1-x_1)x_2,(y_1-x_1)x_3).$$
If we do this and let
$$q(x_1,\psi;\tau)=\tau^{-1}p(x_1,\psi;\tau)=-\frac1{12}\rho(x_1,\psi)\tau^2-\frac1{24}\rho'_{x_1}(x_1,\psi)\tau^3,$$
then, in order to show that $\Pi^0_r$ has at most a folding singularity at a point $(x,y_1,\xi)$, it suffices to show that, at the corresponding point $z$, 
$\text{Hess }\kappa_r(z)\ne 0$ if
\begin{multline*}
\kappa_r(z)=\bigl(z_2(1+r_3/(y_1-z_1))-z_3(q+r_5/(y_1-z_1)),z_3(1+r_3/(y_1-z_1))+z_2(q+r_5/(y_1-z_1)),
\\
(y_1-z_1)^{-1}[(1+\frac{\partial r_3}{\partial y_1})(z_2,z_3)\cdot\xi+(\frac{\partial p}{\partial y_1}+\frac{\partial r_5}{\partial y_1})(-z_3,z_2)\cdot\xi]\bigr).
\end{multline*}
Here, for brevity, $q$ and $p$ denote the functions evaluated at $(z_1,\psi;y_1-z_1)$.

The other projection can be handled in a similar manner.  If we let
$$\eta=(y_1-x_1)\xi,$$
then $\Pi^0_l$ will have at most a folding singularity at $(x,y_1,\xi)$ if the map
\begin{align*}
\kappa_l(\eta,y_1)
 =\Bigl( \, &(1+r_3/(y_1-x_1)) \eta_1+ (q+r_5/(y_1-x_1)) \eta_2, 
\\
&(1+r_3/(y_1-x_1)) \eta_2-(q+r_5/(y_1-x_1)) \eta_1,
\\
& (x_1-y_1)^{-1}[(\eta_1,\eta_2)\cdot(x_2,x_3)-\frac{\partial p}{\partial x_1} (-\eta_2,\eta_1)\cdot(x_2,x_3)]\Bigr)
\\
&\qquad +(x_1-y_1)^{-1} (\frac{\partial r_3}{\partial x_2},\frac{\partial r_3}{\partial x_3}, \frac{\partial r_3}{\partial x_1})(x_2,x_3)\cdot\eta 
\\
&\qquad +
(x_1-y_1)^{-1} (\frac{\partial r_5}{\partial x_2},\frac{\partial r_5}{\partial x_3}, \frac{\partial r_5}{\partial x_1})(-x_3,x_2)\cdot\eta
\end{align*}
has at most folding singularities at the corresponding point $\eta$.

The desired calculations for the main terms of the projections are then summarized in the following

\ble\label{lemma3.7}  The Jacobians of $\kappa_r$ and $\kappa_l$ always have rank at least $2$ if $0\ne x_1-y_1=z_1-y_1$ is small.  Moreover, if $\kappa_r'$ is singular at $z$ and if $|\xi|=1$ then
$$\text{Hess }\kappa_r(z)=|(z_2,z_3)|\, \bigl|\frac{\partial^2q}{\partial z_1^2} - \frac{\partial^3p}{\partial z_1^2\partial y_1}\bigr| +O(|(z_2,z_3)|\rho(y_1-z_1))+O(|(z_2,z_3)|(y_1-z_1)^2),$$
while if $\kappa'_l$ is singular at $(\eta,y_1)$ and $|\eta|=1$ then
$$\text{Hess }\kappa_l(\eta,y_1)=
\bigl|\frac{\partial^2q}{\partial y_1^2}-\frac{\partial^3p}{\partial x_1\partial y^2_1}\bigr| +O(\rho(y_1-z_1))+O((y_1-z_1)^2),$$
assuming in both cases that $0\ne x_1-y_1=z_1-y_1$ is small and that $|(x_2,x_3)|=|(z_2,z_3)|/|y_1-z_1|$ is bounded away from zero.
\ele

\noindent{\bf Remark.}  A straightforward calculation shows that
\begin{equation}\label{3.39}
\frac{\partial^2q}{\partial x_1^2}-\frac{\partial^3p}{\partial x_1^2\partial y_1}=\frac13\rho+\frac1{12}\rho'_{x_1}\cdot(y_1-x_1)
\end{equation}
and
\begin{equation}\label{3.40}
\frac{\partial^2q}{\partial y_1^2}-\frac{\partial^3p}{\partial x_1\partial y_1^2}=-\frac23\rho-\frac3{4}\rho'_{x_1}\cdot(y_1-x_1).
\end{equation}
Our chaotic curvature assumption that $|\rho|+|\rho'|\ne0$ implies that one of these two must be nonzero at a given point if $0\ne x_1-y_1$ is small.  Thus, Lemma \ref{lemma3.7} implies that if we stay near $\gamma_0$ then at points where $\Pi^0_r$ and $\Pi^0_l$ are singular at least one of the projections must have a folding singularity.

\medskip

Let us now present the somewhat tedious proof of this lemma.  Fortunately since the main part of $\Xi$, $\Xi_0$, is linear in $x_2$ and $x_3$, the calculations needed for the scaled geodesics in \eqref{3.25} will also follow from this model case.  We shall say more about this after the proof is complete.

\medskip

\noindent{\it Proof of Lemma \ref{3.7}.}  Let us first handle $\kappa_r$.  To study its Jacobian we first note that in our coordinates
$$|(y_1-z_1)^2D^n_{y_1,z}z_jr_3/(y_1-z_1)|+|D^n_{y_1,z}z_jr_5/(y_1-z_1)|\le C(y_1-z_1)^{5-|n|}, \, \, \,
|n|\le 2, \, \, j=2,3,$$
since $|z_2|,|z_3|\approx |y_1-z_1|$ in view of our assumption that $|(x_2,x_3)|=|\tan\theta|$ is bounded away from zero.  Also,
$$|(y_1-z_1)D^n_{y_1,z}q|+|D^n_{y_1,z}p|\le C|y_1-z_1|^{3-|n|}, \quad |n|\le 2,$$
while if the derivatives just involve $y_1$ and $z_1$ we get an improvement if $\rho$ is small:
$$|(y_1-z_1)D^n_{y_1,z_1}q|+|D^n_{y_1,z_1}p|\le C(|\rho|\, |y_1-z_1|^{3-|n|}+|y_1-z_1|^{4-|n|}), \quad |n|\le 2.$$

With this in mind one checks that the first two columns of $\kappa_r'$ are of the form
$$
\left(
\begin{matrix}
-z_3(\frac{\partial q}{\partial z_1}+O((y_1-z_1)^3))+z_2\frac{\partial}{\partial z_1}r_3/(y_1-z_1) &1+O((y_1-z_1)^2)&O((y_1-z_1)^2)
\\
z_2(\frac{\partial q}{\partial z_1}+O((y_1-z_1)^3))
 +z_3\frac\partial{\partial z_1}r_3/(y_1-z_1)&O((y_1-z_1)^2) &1+O((y_1-z_1)^2)
\end{matrix}
\right),
$$
while its $3,1$-component is of the form
\begin{multline*}
(y_1-z_1)^{-2}\bigl[(z_2,z_3)\cdot\xi +\frac{\partial p}{\partial y_1}(-z_3,z_2)\cdot\xi\bigr]+(y_1-z_1)^{-1}\frac{\partial^2p}{\partial z_1\partial y_1}(-z_3,z_2)\cdot\xi
\\
+O((z_2,z_3)\cdot\xi)+O((y_1-z_1)^2)(-z_2,z_3)\cdot\xi),
\end{multline*}
and since $|(z_2,z_3)|\approx |y_1-z_1|$ its $3,2$ and $3,3$ components are
$$(y_1-z_1)^{-1}[\xi_1+\frac{\partial p}{\partial y_1}\xi_2]+O((y_1-z_1)^2|\xi|)$$
and
$$(y_1-z_1)^{-1}[\xi_2-\frac{\partial p}{\partial y_1}\xi_1]+O((y_1-z_1)^2|\xi|),$$
respectively.  Based on this, if we assume as we may that $z_3=0$, then the determinant of $\kappa_r'$ must be of the form
\begin{multline*}
(1+O((y_1-z_1)^2))\bigl[(y_1-z_1)^{-2}\bigl(
z_2\xi_1
+\frac{\partial p}{\partial y_1}z_2\xi_2+(y_1-z_1)\frac{\partial^2p}{\partial y_1\partial z_1}z_2\xi_2\bigr)\bigr]
\\
-(1+O((y_1-z_1)^2))z_2\frac{\partial q}{\partial z_1}(y_1-z_1)^{-1}\xi_1 +O(\rho z_2|\xi|\, (y_1-z_1))
+O(z_2|\xi|\, (y_1-z_1)^2).
\end{multline*}
From this we deduce that $\text{det }\kappa_r'\ne0$ unless $|z_2\xi_1|$
is smaller than a fixed multiple of $(y_1-z_1)^2|z_2\xi|$ if $(y_1-z_1)$ is small.  Let us therefore assume that $\xi_2=1$.  A more precise calculation then gives that
\begin{multline}\label{3.42}
\xi_1=(y_1-z_1)\bigl(\frac{\partial q}{\partial z_1}-(y_1-z_1)^{-1}\frac{\partial p}{\partial y_1}-\frac{\partial^2p}{\partial y_1\partial z_1}\bigr)+O(\rho(y_1-z_1)^3)+O(y_1-z_1)^4)
\\
\text{if } \, \, \text{det }\kappa_r'=0 \, \, \text{and } \, \xi_2=1.
\end{multline}
Clearly, $\text{rank }\kappa_r'\ge2$ everywhere if $(y_1-z_1)$ is small.

Assuming \eqref{3.42}, let us compute $\text{Ker }\kappa_r'$.  If $X$ is a unit vector in $\text{Ker }\kappa'_r$ then clearly if $z_3=0$ its second and third components must be $O((y_1-z_1)^2)$ and $O(y_1-z_1)$, respectively.  More precisely, if we assume that the first component is $-1$, then
$$X=\bigl(
-1,\, O((y_1-z_1)^2) ,\, z_2\frac{\partial q}{\partial z_1}+O(\rho(y_1-z_1)^3)+O((y_1-z_1)^4)
\bigr).
$$
Consequently,
\begin{multline*}
\bigl(<X,\nabla>\bigr)^2=\frac{\partial^2}{\partial z_1^2}+\bigl(1+O(\rho(y_1-z_1))+O((y_1-z_1)^2)\bigr)\times
2z_2\frac{\partial q}{\partial z_1}\frac{\partial^2}{\partial z_1\partial z_3} 
\\
+O((y_1-z_1)^2)\frac{\partial^2}{\partial z_1\partial z_2}+O((y_1-z_1)^4)\nabla^2.
\end{multline*}

To compute the cokernel we note that the last two rows of the transpose of the Jacobian are of the form
$$
\left(
\begin{matrix}
1+O((y_1-z_1)^2)&O((y_1-z_1)^2)&(y_1-z_1)^{-1}(\xi_1+\xi_2\frac{\partial p}{\partial y_1})+O((y_1-z_1)^2)
\\
O((y_1-z_1)^2)&1+O((y_1-z_1)^2) &(y_1-z_1)^{-1}(\xi_2-\xi_1\frac{\partial p}{\partial y_1})+O((y_1-z_1)^2)
\end{matrix}
\right).$$
Based on this, if $Y$ is a unit vector in the cokernel, then its first component must be $O(\rho(y_1-z_1))+O((y_1-z_1)^2)$.  Since we are assuming that $\xi_2=1$ we can say more using the last row of the transpose.  Namely, a vector of the form
$$Y=\bigl(O(\rho(y_1-z_1))+O((y_1-z_1)^2), \, -1,\, y_1-z_1+O((y_1-z_1)^2)\bigr)$$
is in the cokernel.

Let us now evaluate $<Y,\, (<\nabla,X>)^2\kappa_r>$ to compute the Hessian.  We first compute the contribution to the Hessian of the third component when, as above, $\xi_2=1$
\begin{multline*}
(<\nabla,X>)^2\kappa_r^3
=(1+O(\rho(y_1-z_1))+O((y_1-z_1)^2))\times
\\
\Bigl[\frac{-2z_2}{(y_1-z_1)^2}\frac{\partial q}{\partial z_1}+\frac{z_2}{(y_1-z_1)^2}\frac{\partial^3p}{\partial y_1\partial z_1^2}
+\frac2{(y_1-z_1)^3}(z_2\xi_1+z_2\frac{\partial p}{\partial y_1})
\\
+\frac{2z_2}{(y_1-z_1)^2}\frac{\partial^2p}{\partial y_1\partial z_1}\Bigr]
+O(z_2(y_1-z_1)).
\end{multline*}
If we recall \eqref{3.42}, we conclude that
$$(<\nabla,X>)^2\kappa^3_r=\frac{\partial^3p}{\partial y_1\partial z_1^2}\frac{z_2}{y_1-z_1}+O(\rho(y_1-z_1))+O((y_1-z_1)^2).$$
Similar considerations yield
$$\bigl(<\nabla,X>\bigr)^2\kappa^2_r=z_2\frac{\partial^2q}{\partial z_1^2}+O(\rho(y_1-z_1)^2)+O((y_1-z_1)^3),$$
and
$$\bigl(<\nabla,X>\bigr)^2\kappa^1_r=O((y_1-z_1)).$$

If we recall the form of $Y$, we conclude that
$$<Y,(<\nabla,X>)^2\kappa_r>=z_2\bigl(\frac{\partial^3p}{\partial y_1\partial z_1^2}-\frac{\partial^2q}{\partial z_1^2}\bigr)+O(\rho((y_1-z_1)^2)+O((y_1-z_1)^3),$$
which gives rise to the first part of Lemma \ref{lemma3.7} since $|Y|=1+O(y_1-z_1)$.

\medskip

To handle the second half of Lemma \ref{3.7} we should first notice that the main part of the map $\kappa_l$ is basically the same as that of $\kappa_r$ with the roles of $(x_2,x_3)$ and $\eta$ and $x_1$ and $y_1$ reversed.  Thus it should not be surprising that the preceding arguments allow us to compute $\text{Hess }\kappa_l$ when $\kappa_l'$ is singular.

Indeed, if we argue as before, we find that when $\eta=(0,1)$ the determinant of $\kappa_l'$ equals
\begin{multline*}
(x_1-y_1)^{-2}(x_3+x_2\frac{\partial p}{\partial x_1})+(x_1-y_1)^{-1}\frac{\partial^2p}{\partial x_1\partial y_1}x_2-(x_1-y_1)^{-1}x_2\frac{\partial q}{\partial y_1}
\\
+O(\rho(x_1-y_1))+O((x_1-y_1)^2),
\end{multline*}
with $x'=(x_2,x_3)$.  Consequently,
$$
x_3=\frac{\partial q}{\partial y_1}x_2(x_1-y_1)-\frac{\partial^2 p}{\partial x_1\partial y_1}x_2(x_1-y_1)-x_2\frac{\partial p}{\partial x_1}+O(\rho(x_1-y_1)^3)+O((x_1-y_1)^4),
$$
when $\text{det }\kappa'_l=0$ and $\eta=(0,1)$.
Furthermore, under these assumptions
$$X=\bigl(-\frac{\partial q}{\partial y_1}+O((x_1-y_1)^3),\, O(x_1-y_1), \, 1\bigr)
\in \text{Ker }\kappa'_l,
\quad \kappa'_l=\partial\kappa_l/\partial(\eta,y_1),
$$
and
$$Y=\bigl(-x_2,\, O(\rho(x_1-y_1))+O((x_1-y_1)^2), \, (x_1-y_1)+O((x_1-y_1)^2)\bigr)\in \text{Ker }(\kappa_l')^t.$$
Since then
$$<Y,\, (\nabla,X>)^2\kappa_l>=x_2(\frac{\partial^3 p}{\partial x_1\partial y_1^2}-\frac{\partial^2 q}{\partial y_1^2}) +O(\rho(x_1-y_1))+O((x_1-y_1)^2),$$
we get the second part of the lemma after noting that $|Y|=|x_2|+O(x_1-y_1)$.  

This completes the proof of Lemma \ref{3.7}.

\noindent{{\it Proof of Proposition \ref{prop3.6}}  We need to check that we can construct weights $a_{\gamma_{x',s}}$ so that \eqref{3.22} holds and so that the scaled weights \eqref{3.26} belong to a bounded subset of $C^\infty$ and moreover $\Pi^k_r$ and $\Pi^k_l$ have at most folding singularities at points where $a_{\gamma^k_{x',s}}(y)\ne 0$ and $y\in \gamma^k_{x',s}$.  Here $\gamma^k_{x',s}$ is the scaled geodesic as in \eqref{3.25}, while $\Pi^k_r$ and $\Pi^k_l$ are the associated right and left projections of the associated canonical relation ${\cal C}_k$.

If as above we parameterize the unscaled geodesics by variables $(x_1,\theta,\psi)$ reflecting the intersection with $\gamma_0$ and the resulting initial unit tangent vector we then as before let $(x_2,x_3)=\tan\theta(\cos\psi,\sin\psi)$.  It then follows that $\Pi^k_r$ and $\Pi^k_l$ must be the analogs of the projections $\Pi^0_r$ and $\Pi^0_l$ just studied, where $\Pi_0$ is replaced by
$$\Xi_k(x,y_1)=2^k\Xi(x_1,2^{-k}x_2,2^{-k}x_3,y_1).$$
If we recall \eqref{3.3.38} and note that $\Xi_0$ is linear in $(x_2,x_3)$, we conclude that
$$\Xi_k=\Xi_0+O(2^{-k}(x_1-y_1)^3),$$
where $2^k$ times the error term belongs to a bounded subset of $C^\infty$.  Hence, if replace $r_3$ and $r_5$ in the preceding arguments by error terms of the form $r_3+O(2^{-k}(x_1-y_1)^3)$ and $r_5+O(2^{-k}(x_1-y_1)^3)$, we can argue as above to see that for the resulting analogs $\kappa^k_r$ and $\kappa^k_l$ of $\kappa_r$ and $\kappa_l$ we have
\begin{multline}\label{3.43}
\text{Hess }\kappa^k_r=|(z_2,z_3)|\, 
\bigl|\frac{\partial^2q}{\partial x_1^2}- \frac{\partial^3p}{\partial x_1^2\partial y_1}\bigr| 
\\
+|(z_2,z_3)|\bigl(O(\rho(x_1-y_1))+O((y_1-x_1)^2)+O(2^{-k})\bigr),
\end{multline}
if $(\kappa^k_r)'$ is singular at $z=(x_1,z_2,z_3)$ and $|\xi|=1$, as well as
\begin{equation}\label{3.44}
\text{Hess }\kappa^k_l=
\bigl|\frac{\partial^2q}{\partial y_1^2}-\frac{\partial^3p}{\partial x_1\partial y_1^2}\bigr| 
+O(\rho(y_1-x_1))+O((y_1-x_1)^2)+O(2^{-k}),
\end{equation}
if $(\kappa^k_l)'$ is singular at $(\eta,y_1)$ and $|\eta|=1$.

To proceed, we recall that we may assume that $2^{-k}$ is as small as we wish.  This corresponds to making the parameter $r$ in Proposition \ref{prop3.6} small.  We also need to recall that our variable curvature assumption amounts to the condition that for some $c_0>0$ the coefficients of $p$ and $q$ satisfy
$$c_0\le |\rho|+|\rho'_{x_1}|\le c_0^{-1}.$$
Recalling \eqref{3.39} and \eqref{3.40}, if as above $\alpha$ denotes the lengths of our geodesics, let us choose $\alpha_0\le \alpha/2$ so that the ``quadratic'' error terms in \eqref{3.43} and \eqref{3.44} both satisfy
$$|O(y_1-x_1)^2)|\le c_0/10\quad \text{if } \, \, |x_1-y_1|\le \alpha_0.$$
If we then also assume that $k$ is large enough so that the $O(2^{-k})$ error terms satisfy $O(2^{-k})\le c_0/10$, we conclude that we can choose $c_1>0$ and $0<\alpha_1<\alpha_0$ so that
\begin{multline*} \bigl|\frac{\partial^2q}{\partial x_1^2}- \frac{\partial^3p}{\partial x_1^2\partial y_1}\bigr|, \, \bigl|\frac{\partial^2q}{\partial y_1^2}-\frac{\partial^3p}{\partial x_1\partial y_1^2}\bigr| \, \ge \frac{c_0}{10}|x_1-y_1|, 
\\
\text{if } \, \, |\rho|\le c_1\, \, \, \text{and } \, |x_1-y_1|\in [\alpha_1,\alpha_0].
\end{multline*}
For the remaining case, we need to assume further that the $2^{-k}$ error terms satisfy $|O(2^{-k})|\le c_1/10$.  In this case there must be a constant $\alpha_2\le \alpha_1$ so that
$$ \bigl|\frac{\partial^2q}{\partial x_1^2}- \frac{\partial^3p}{\partial x_1^2\partial y_1}\bigr|, \, \bigl|\frac{\partial^2q}{\partial y_1^2}-\frac{\partial^3p}{\partial x_1\partial y_1^2}\bigr| \, \ge c_1/10, 
\quad
\text{if } \, \, |\rho|\ge c_1, \, \, \text{and }\, |x_1-y_1|\le \alpha_2$$

Because of this we can clearly choose weights with the required properties.  For the first case we choose a bump function $\beta_1\in C^\infty_0((\alpha_1,\alpha_0))$ which equals one in the middle half of $(\alpha_1,\alpha_0$, while for the second case we choose $\beta_2\in C^\infty_0((0,\alpha_2))$ satisfying $\beta_2=1$ on $[\alpha_2/4,\alpha_2/2]$.  If we then let $a_{\gamma_x}(y)$ be equal to $\beta_2(|x_1-y_1|)$ if $|\rho(x_1)|\ge c_1$ and equal to $\beta_1(|x_1-y_1|)$ otherwise, it then follows that the resulting pullback $a_{\gamma_{x',s}}$ to the $(x',s)$ coordinates will have the desired properties.  This finishes our proofs.

\newsection{More general lower bounds in variably curved manifolds}

In this section we shall briefly indicate how we can obtain slightly less favorable lower bounds on the dimension of a larger collection of sets.  Specifically we have the following analog of Theorem 3.1.

\bth\label{thm4.1}  Assume that $M^3$ is variably curved in the sense of Definition 3.1.  Then
\begin{equation}\label{4.1}  \text{dim }\Omega\ge 9/4
\end{equation}
if $\Omega\subset M^3$ satisfies $|\cup_{\lambda>0}\Omega^*_{\alpha,\lambda}|>0$ for a given sufficiently small $\alpha$ with $\Omega^*_{\alpha,\lambda}$ being as in (1.1).
More generally, if ${\cal F}$ is a family of geodesics of length $\alpha$ satisfying the chaotic curvature condition then the same conclusion holds if
$$|\{x\in M^3: \, \exists \gamma_x\in {\cal F} \, \, \text{with } \, \, x\in \gamma_x \, \, \text{and } \, \, |\gamma_x\cap\Omega|>0\}|>0.
$$
\eth

For the proof one first notices that the result follows from certain estimates for the Nikodym maximal operator in \eqref{nikmax2}.  In this case, \eqref{4.1} would follow if we could show that if $\varepsilon>0$ is fixed then for every $0<\lambda<1$ there is a constant $C_{\varepsilon,\lambda}$ so that
\begin{equation}\label{4.2}
\delta^{3/4+\varepsilon}|\{x: \, f^*_\delta(x)\ge \lambda\}|\le C_{\varepsilon,\lambda}|E|, \quad f=\chi_E.
\end{equation}
The dependence on $\lambda$ is not important for our application; however, we should point out that the present methods yield much worse bounds than the ones $C_{\lambda,\varepsilon}\le C_\varepsilon \lambda^{-5/2}$ for the constant curvature case (where a more favorable dependence on $\delta$ was also obtained).

To prove \eqref{4.2} we shall appeal to an auxiliary maximal operator which is just a truncated version of the one in (2.9).  Specifically, for a given $0<\lambda<1$, we shall consider
\begin{equation}\label{4.3}
( {\cal A}_{\lambda,\delta}f)(x')=\sup_{\{\gamma_{x'}\ni (0,x'): \gamma_{x'}\cap\gamma_0\ne \emptyset\}} |T^\delta_{\gamma_{x'}}|^{-1} \int_{\{y\in T^\delta_{\gamma_{x'}}: \text{dist }(y,\gamma_0)\ge \lambda\}} |f(y)|\, dy,
\end{equation}
when $\gamma_0$ satisfies the chaotic curvature condition.
As in the statement of the theorem, we assume here that the length $\alpha$ of the tubes $T^\delta_{\gamma_{x'}}$ is small.  By staying away from the common geodesic $\gamma_0$ we can avoid small scale issues which complicate the analysis since near $\gamma_0$ the geometry looks Euclidean.

This truncation also allows us to exploit the remark after Lemma 3.7 which said that the underlying Fourier integral operators that govern the averages in \eqref{4.3} have canonical relations with at most one-sided folding singularities if $\alpha$ is small.  Greenleaf and Seeger \cite{GS} showed that such Fourier integral operators are bounded on $L^2$ with a loss of $1/4$ derivative.  By appealing to this result we can argue as in the last section to conclude that if $r>0$ is small
\begin{equation}\label{4.4}
\|{\cal A}_{\lambda,\delta}f\|_{L^2(|x'|<r)}\le C_{\lambda,\varepsilon}\delta^{-1/4-\varepsilon}\|f\|_2.
\end{equation}

From this we can obtain \eqref{4.2} using our earlier arguments.  To see this, we first recall that estimates like \eqref{4.4} are used to obtain lower bounds for $|E|$ when we assume that there are many tubes $T^\delta_j$ intersecting $\gamma_0$ for which
\begin{equation}\label{4.5}
|E\cap T^\delta_j|\ge \lambda |T^\delta_j|.
\end{equation}
If we replace $\lambda$ by $c\alpha\lambda$ in \eqref{4.4} then we could use the resulting inequality to obtain the desired lower bounds for $|E|$ if we knew, say, that
$$|\{y\in E: \, \text{dist }(y,\gamma_0)\ge c\alpha\lambda\}\cap T^\delta_j|\ge \frac\lambda2 |T^\delta_j|.
$$
But this of course follows from \eqref{4.5} since 
$$|\{y\in T^\delta_j: \, \text{dist }(y,\gamma_0)\le c\alpha\lambda\}|\le \frac\lambda2 |T^\delta_j|$$
if $c$ is small.

\newsection{Negative results and some problems}

Let us begin this section by showing how the maximal estimates (2.3) and (2.11) for the constant curvature case break down if one is working in a variably curved manifold.  The arguments are similar to those given for more degenerate situations by Minicozzi and the author \cite{MS}.

To provide counterexamples we shall fix a common geodesic $\gamma_0$ and work in Fermi normal coordinates about it.  Recall then that given any $x_1$ we can find a $\psi=\psi(x_1)$ so that $\rho(x_1,\psi)=0$ if, as in \eqref{3.14}, $\rho$ is $3!$ times the main Taylor coefficient of the component of $\gamma_{x_1\theta\psi}$ which is orthogonal to the Fermi two plane of tangency at $\gamma_{x_1\theta\psi}\cap\gamma_0$.  If we fix $x_1=\overline{x}_1$, then after perhaps rotating our coordinates around $\gamma_0$, we may assume for simplicity that $\psi=0$.  

To proceed, as in Proposition 3.2, let $\gamma_{x_1\theta}(t)$ be the geodesic parameterized by arclength satisfying
$$\gamma_{x_1\theta}=(x_1,0,0) \quad \text{and } \, \, \dot \gamma_{x_1\theta}(0)=(\cos\theta,\sin\theta,0).$$
We then set
$\kappa(x_1,\theta,t)=\gamma_{x_1\theta}(t).$
Using Proposition 3.2, one can check that there are $\delta_j>0$ so that the Jacobian of $\kappa$ satisfies
\begin{equation}\label{5.1}
|\text{det }\kappa'|\approx \theta \quad \text{if } \, \,
|x_1-\overline{x}_1|+|\theta|\le \delta_1 , \, \, \, \delta_2\le |t|\le 2\delta_2.
\end{equation}

Based on this, if $f^*_\delta$ is as in \eqref{nikmax2} we can easily show that (2.3) breaks down if the family of geodesics involved satisfies the chaotic curvature condition (see Definition 3.1).  Such an example would be where ${\cal F}$ is a family of geodesics which are close to $\gamma_0=\{(x_1,0,0): \, 0\le x_1\le \alpha\}$ if one considers the metric
\begin{equation}\label{5.2}
dx^2+\bigl((x^2-x^2_3)\cos x_1 + 2 x_2x_3\sin x_1 \bigr) dx^2_1
\end{equation}
on ${\Bbb R}^3$.  For this example the $O(\theta^2)$ ``error'' terms in \eqref{3.12} vanish for $\gamma_0$ allowing worse counterexamples.  

To be more specific, if we consider \eqref{5.2}, let us take $\overline{x}_1=0$ here since $g_{11,23}(0,0,0)=0$.  We then let $f_\delta=\chi_{\Omega_\delta}$, where
$$\Omega_\delta =\{x: \, |x_1|+|x_2|\le \delta^{1/4}, \, \, |x_3|\le \delta\}.$$
Then clearly
$$\|f\|_p\approx \delta^{3/2p}.$$
On the other hand, using Proposition 3.2 and \eqref{5.1} with, say $\delta_1/2\le|\theta|\le \delta_1$, one sees that there must be a set $\Omega^*_\delta$ of measure larger than a fixed multiple of $\delta^{1/4}$ so that $f^*_\delta(x)\ge c\delta^{1/4}$ if $x\in \Omega^*_\delta$, if, as in (2.1), $f^*_\delta$ is the Nikodym maximal function of $f=f_\delta$.  Consequently, 
$$\|f_\delta^*\|_q\ge c\delta^{1/4+1/4q}.$$
Based on this, we immediately see that (2.3) cannot hold since when $p=5/2$ and $q=10/3$, 
$$\|f^*_\delta\|_{10/3}\, /\, \|f_\delta\|_{5/2}\ge c\delta^{-11/40},$$
while the bounds (2.3) for the constant curvature case say that the ratio is $O(\delta^{-1/5-\varepsilon})$ for any fixed $\varepsilon>0$.  Curiously, if one weakens (2.3) by replacing the norm in the left by $L^p$, bounds like those obtained by Bourgain \cite{B1}
$$\|f^*_\delta\|_p\le C_\varepsilon \delta^{1-3/p-\varepsilon}\|f\|_p, \quad 1\le p\le 7/3$$
would be best possible under the present hypotheses.

The assumption that a given geodesic should satisfy the chaotic curvature condition 
is all that was used in the proof of all of the results for the variably curved case.  Since we avoided small scales $\lambda$ in \eqref{3.16} we were able to essentially ignore the $O(\theta^2)$ terms in \eqref{3.12} and \eqref{3.13}.  The next counterexample along with the preceding suggests that it might be necessary to use these terms to improve the results of the last section.

To be more specific, let us consider a general variably curved manifold as in Definition 3.1.  We shall assume as above that $g_{11,23}=0$ at $(\overline{x}_1,0,0)$.  However, since we are not assuming now that $g_{12,23}$ vanishes there, in view of \eqref{3.5}, we must modify the above counterexample.  Assuming as we are that $g_{12,231}\ne 0$ at this point we must only consider very small values of $\theta$ so that the fourth order terms dominate the third order terms in the Taylor expansion of the last coordinate of $\gamma_{x_1\theta}$.  To achieve this, we now let $f_\delta=\chi_{\Omega_\delta}$, where
$$\Omega_\delta=\{x: \, |x_1-\overline{x}_1|+|x_2|\le \delta^{1/5}, \, \, |x_3|\le 2\delta \}.$$
One can then see from Proposition 3.2 that there must be a constant $c>0$ so that 
$$\gamma_{x_1\theta}(t)\in \Omega_\delta\quad \text{if } \, \, |x_1-\overline{x}_1|\le c\delta^{1/5} \, \, \, \text{and } \, \, |\theta|\le c\delta^{2/5}.$$
Based on this \eqref{5.1} ensures that
$$f^*_\delta(x)\ge c'\delta^{1/5} \quad\text{if } \, \, x\in \Omega^*_\delta
\, \, \, \text{with } \, \, |\Omega^*_\delta|\ge c'\delta^{3/5}.$$
Consequently,
$$\|f^*_\delta\|_q\, /\, \|f\|_p \ge c\delta^{1/5+3/5q-8/5p}.$$
Because of this, one immediately sees that (2.3) cannot hold in any variably curved manifold.  Indeed, if for $p=5/2$ one wishes to have the Nikodym maximal operator bounded from $L^p$ to $L^q$ with norm $O(\delta^{-1/5-\varepsilon})$ for all $\varepsilon>0$, then one must take $q=p=5/2$.  Moreover, if the maximal operator is bounded from $L^p$ to $L^p$ with norm $O(\delta^{1-3/p-\varepsilon})$ for all $\varepsilon>0$, then $p\le 5/2$.

One can also use this construction of course to show that the bounds (2.11) for the auxiliary maximal operator defined in (2.9) cannot hold in variably curved manifolds.  Moreover, if one wishes for bounds like (2.11) to hold here, then one must modify the definition (2.9) replacing the weights in (2.10) by
$$\bigl(\text{dist }(y,\gamma_{x'}\cap\gamma_0)\bigr)^{3/2}.$$
This would be the smallest power of this distance function for which bounds like (2.11) could hold.  If one could prove the inequality for this mollified auxiliary operator, then the proof of (2.3) would give the bounds
$$\|f^*_\delta\|_{14/3}\le C_\varepsilon \delta^{-1/7-\varepsilon}\|f\|_{7/2},$$
and using this one would see that (2.5) must also hold when $M^3$ is variably curved.

It would also be interesting to study what happens in higher dimensions.  Here it would certainly be of interest to improve on the lower bound $\text{dim }\Omega \ge (n+1)/2$ for Nikodym-type subsets of $n$-dimensional subsets of symmetric spaces (see \cite{MS}).  The results of Section 2 cover the case of $3$-dimensional symmetric spaces since such a manifold must have constant curvature.  This of course does not happen in higher dimensions.  A typical example is $ {\Bbb C}P^n$ where the curvature is not constant and for related reasons there are some totally geodesic submanifolds but not nearly as in Euclidean spaces of the same dimension.  Thus, ${\Bbb C}P^n$ would in some sense represent an intermediate case between the types of manifolds considered in Sections 2 and 3.  For this reason it would be interesting to see whether the analog of Wolff's Euclidean bounds in \cite{W} hold.  That is, does one always have that $\text{dim }\Omega\ge (d+2)/2$ if $\Omega\subset {\Bbb C}P^n$ is a Nikodym-type set and $d=2n=\text{dim}_{\Bbb R}{\Bbb C}P^n$?  Along the same lines, does this lower bound always hold if $M^d$ is an Einstein manifold of dimension $d>3$?  It would also be interesting to try to formulate a condition in the spirit of Section 3 which would ensure that one is far from the symmetric case and also allow improvements over the easy bounds $\text{dim }\Omega\ge (n+1)/2$ for Nikodym-type subsets.  As was pointed out in \cite{MS}, if $n$ is odd there are always examples where the lower bound cannot be improved so such a condition, like the one in Section 3, would have to rule out these degenerate and hopefully atypical cases.

\end{document}